\documentclass{article}
\usepackage{amsmath}
\usepackage{amsfonts}
\usepackage{csquotes}

\usepackage{amsthm}
\usepackage{amssymb}

\usepackage[citestyle=authoryear,bibstyle=authoryear]{biblatex}
\newtheorem{theorem}{Theorem}[section]
\newtheorem{corollary}{Corollary}[theorem]

\newtheorem{definition}[theorem]{Definition}

\newtheorem{conjecture}[theorem]{Conjecture}
\newtheorem{question}[theorem]{Question}

\newtheorem{remark}[theorem]{Remark}

\title{On the Hierarchy of Natural Theories}

\author{James Walsh
}

\begin{document}

\maketitle

\begin{abstract}
It is a well-known empirical phenomenon that natural axiomatic theories are pre-well-ordered by consistency strength. Without a precise mathematical definition of ``natural,'' it is unclear how to study this phenomenon mathematically. We will discuss the significance of this problem and survey some strategies that have recently been developed for addressing it. These strategies emphasize the role of reflection principles and ordinal analysis and draw on analogies with research in recursion theory.
\end{abstract}


\section{The Consistency Strength Hierarchy}

The starting points of this story are G\"{o}del's incompleteness theorems. These theorems establish fundamental restrictions on what can be proven on the basis of any reasonable axiomatic theory. By a \emph{reasonable} axiomatic theory I mean a consistent, recursively axiomatized theory that interprets a modicum of arithmetic.

\begin{theorem}[\cite{godel1931formal, rosser1936extensions}]
No reasonable axiomatic theory is complete.\footnote{This is technically Rosser's strengthening of G\"{o}del's Theorem, if we demand only that reasonable theories are consistent. G\"{o}del used the assumption of $\omega$-consistency in his proof.}
\end{theorem}

\begin{theorem}[\cite{godel1931formal}]
No reasonable axiomatic theory proves its own consistency.
\end{theorem}

The first incompleteness theorem shows that, no matter what principles we endorse, there will be questions that cannot be resolved on the basis of those principles. This means that there is no universal axiom system within which mathematics can be developed. Instead, we are left with a vast array of axiomatic theories. The second incompleteness theorem yields the primary method for mapping out this vast array of theories. Fix some reasonable base theory $B$. Let's now recall some standard definitions.
\begin{definition}
    For reasonable theories $T$ and $U$, we say that $T \leq_{\mathsf{Con}}U $ if $B$ proves that the consistency of $U$ implies the consistency of $T$. $T\equiv_{\mathsf{Con}}U$ if  $T \leq_{\mathsf{Con}} U $ and $U \leq_{\mathsf{Con}} T$. $T <_{\mathsf{Con}} U$ if $T \leq_{\mathsf{Con}} U$ and $U \nleq_{\mathsf{Con}} T$.
\end{definition}
The notion of strength that $<_{\mathsf{Con}}$ engenders is known as \emph{consistency strength}. The structure of \emph{all} reasonable theories ordered by consistency strength is quite complicated. As far as I know, the following theorems are folklore:

\begin{theorem}[Folklore]
The ordering $<_{\mathsf{Con}}$ is not \emph{pre-linear}, i.e., there are theories $T$ and $U$ such that $T\nleq_{\mathsf{Con}}U$ and $U \nleq_{\mathsf{Con}} T$.
\end{theorem}

\begin{theorem}[Folklore]
The ordering $<_{\mathsf{Con}}$ is \emph{ill-founded}, i.e., there is a sequence $T_0 >_{\mathsf{Con}} T_1 >_{\mathsf{Con}} T_2 >_{\mathsf{Con}} \dots$ where each $T_i$ is consistent.
\end{theorem}

All known instances of non-linearity and ill-foundedness have been discovered by defining theories in an ad-hoc manner using self-reference and other logical tricks. The restriction of $<_{\mathsf{Con}}$ to \emph{natural} axiomatic theories---those that arise in practice---is a \emph{pre-well-ordering}.\footnote{This is to say that the induced ordering on the equivalence classes modulo equiconsistency is a well-ordering.} Here is a tiny snapshot of such theories:
$$\mathsf{EA}, \mathsf{EA}^+, I\Sigma_n, \mathsf{PA}, \mathsf{ATR}_0, \Pi^1_n\mathsf{CA}_0, \mathsf{PA}_n, \mathsf{ZF}, \mathsf{AD}^{L(\mathbb{R})}$$
The well-ordering phenomenon persists, taking a very liberal view of what constitutes a ``natural'' theory. Note that the theories just cited run the gamut from weak fragments of arithmetic to subsystems of analysis and all the way to strong extensions of set theory. These theories come from different areas of mathematics (e.g., arithmetic, analysis, set theory) and often codify different conceptions of mathematics. Indeed, many of the natural extensions of set theory that have been investigated are jointly inconsistent, yet comparable according to consistency strength.

Explaining the contrast between natural axiomatic theories and axiomatic theories in general is widely regarded as a major outstanding conceptual problem in mathematical logic. The following passage is representative:
\begin{displayquote}[\cite{friedman2013slow}, p. 382]
The fact that ``natural'' theories, i.e. theories which have something like an ``idea'' to them, are almost always linearly ordered with regard to logical strength has been called one of the great mysteries of the foundations of mathematics.
\end{displayquote}
One can find similar quotes elsewhere in the literature. For instance, Shelah writes that explaining ``the phenomena of linearity of consistency strength'' is a ``logical dream'' \parencite{shelah2002logical}.\footnote{These two quotes emphasize linearity, but other logicians also emphasize well-foundedness \parencite{koellner2010independence}.}

If it is true that natural axiomatic theories are pre-well-ordered by consistency strength, and not an illusion engendered by a paucity of examples, then one would like to prove that it is true. However, the claim that natural axiomatic theories are pre-well-ordered by consistency strength is not a strictly mathematical claim. The problem is that we lack a definition of ``natural axiomatic theory.'' Without a precise definition of ``natural,'' it is not clear how to prove this claim, or even how to state it mathematically.

Here is my plan for the rest of the paper. In \textsection \ref{new-axioms} I will discuss the relevance of the pre-well-ordering phenomenon to the search for new axioms in set theory. In \textsection \ref{section on ordinal analysis}--\ref{Ordinal analysis via iterated reflection} I will discuss two research programs that emerged from the incompleteness theorems, namely, ordinal analysis and the study of iterated reflection principles. Both programs will play an important role in the approaches to the pre-well-ordering problem discussed in this paper. In \textsection \ref{turing-degrees} I will discuss a phenomenon from recursion theory, namely, the well-ordering of natural Turing degrees by Turing reducibility; this phenomenon is analogous in many ways to the pre-well-ordering of natural theories by consistency strength. In \textsection \ref{operators}--\ref{reflection-and-ordinals} I will discuss approaches to the pre-well-ordering problem that have emerged in recent work. Finally, in \textsection \ref{conclusions} I will discuss remaining questions and directions for future research.

\section{The Search for New Axioms}\label{new-axioms}

In this section I will discuss how the pre-well-ordering phenomenon informs the search for new axioms in set theory.

Set theory has been developed in an explicitly axiomatic fashion, at least since Zermelo isolated the axioms that are the precursors to the standard $\mathsf{ZFC}$ axioms. The progenitors of set theory discovered that various mathematical structures could be realized in set-theoretic terms and that their properties could be established on the basis of set-theoretic reasoning. Logicians often formulate this aspect of set theory in terms of the technical notion of \emph{interpretability}. Roughly, a theory $T_1$ is interpreted in another theory $T_2$ if there is a translation $\tau$ from the language of $T_1$ to the language of $T_2$ such that, for each sentence $\varphi$ of the language of $T_1$, if $T_1\vdash\varphi$ then $T_2\vdash \tau(\varphi)$.\footnote{Of course, $\tau$ cannot be an arbitrary function. For a precise definition see \cite[Ch. 7]{lindstrom2017aspects}.} As Koellner notes, interpretability is ubiquitous in mathematics: ``Poincar\'e provided an interpretation of two dimensional hyperbolic geometry in the Euclidean geometry of the unit circle; Dedekind provided an interpretation of analysis in set theory; and G\"odel provided an interpretation of the theory of formal syntax in arithmetic'' \parencite{koellner2010independence}. We may restate the aforementioned discovery that large swathes of mathematics can be developed in set theory as the discovery that $\mathsf{ZFC}$ axioms interpret the axioms of various other branches of mathematics.

Despite the fact that $\mathsf{ZFC}$ interprets many other theories, many of the central problems of set theory (such as the Continuum Hypothesis, the Projective Measure problem, and Suslin's Hypothesis) \emph{cannot} be resolved on the basis of the $\mathsf{ZFC}$ axioms. That is, the $\mathsf{ZFC}$ axioms suffice for resolving mainstream problems in all of the branches of mathematics with the conspicuous exception of set theory itself. This has motivated the search for new axioms for set theory.

Set theorists have investigated many extensions of $\mathsf{ZFC}$, including large cardinal axioms, determinacy axioms, forcing axioms, and more. Is it possible to make rational judgments about these principles and their consequences? Some point to the interpretability strength of $\mathsf{ZFC}$ as its core virtue. Indeed, it is the interpretability strength of $\mathsf{ZFC}$ that allows us to carry out so much mathematics on the basis of set-theoretic principles and to regard $\mathsf{ZFC}$ as a ``foundations'' for mathematics. Perhaps, so the thought goes, we should embrace those axiomatic extensions of set theory that most enhance its interpretability power, thereby enriching the domain of mathematics as much as possible.\footnote{Such a thought echoes Cantor's dictum that mathematicians are totally free in their practice, constrained only by considerations of consistency and definiteness \parencite{cantor1976foundations}.} Steel has promoted the maxim ``maximize interpretability strength'' as a guiding principle in the search for new axioms \parencite{steel2000mathematics}. Let's call this principle \textsc{steel's maxim}.\footnote{In an early draft of this paper I referred to this principle as \textsc{maddy's maxim}, since Maddy promoted a related principles \parencite{maddy1988believing, maddy1988believingII}. However, whereas \textsc{steel's maxim} is a principle to choose \emph{theories} that maximize interpretability strength, the principle Maddy formulated is about maximizing certain objects, namely, sets.} For a large swathe of theories, including extensions of set theory, maximizing interpretability strength coincides with maximizing consistency strength. For a sentence $\varphi$ that is independent of $\mathsf{ZFC}$, one can imagine the following four possibilities:
\begin{enumerate}
\item $\varphi$ increases interpretability strength but $\neg \varphi$ does not.
\item $\neg \varphi$  increases interpretability strength but $\varphi$ does not.
\item Neither $\varphi$ nor $\neg \varphi$ increases interpretability strength.
\item Both $\varphi$ and $\neg \varphi$ increase interpretability strength.
\end{enumerate}

It turns out that all four possibilities are realized.\footnote{Of course, the first is realized if and only if the second is realized, by double-negation elimination.} In the fourth case we \emph{cannot} follow \textsc{steel's maxim} or we would land ourselves in inconsistency. However, this is not typically taken as a refutation of \textsc{steel's maxim}. The key point here is that when one restricts oneself to \emph{natural} theories, only the first three possibilities are realized, since natural theories are linearly ordered by consistency strength. See \cite{koellner2010independence} for a discussion of this point.

Consider once more the variety of extensions of $\mathsf{ZFC}$ that are investigated by set theorists: large cardinal axioms, axioms of definable determinacy, forcing axioms, and more. These axiom systems have different motivations and often codify different foundational conceptions of mathematics. Nevertheless, they are pre-well-ordered by consistency strength. Indeed, they are pre-well-ordered according to all the common notions of proof-theoretic strength (such as 1-consistency strength, $\Pi^1_1$ reflection strength, and so on). A consequence is that these axiom systems all converge on arithmetic statements and even analytic statements. As Steel writes:
\begin{displayquote}[\cite{steel2014godel}, p. 159]
Thus the well-ordering of natural consistency strengths corresponds to a well-ordering by inclusion of theories of the natural numbers. There is no divergence at the arithmetic level, if one climbs the consistency strength hierarchy in any natural way we know of\dots. Natural ways of climbing the consistency strength hierarchy do not diverge in their consequences for the reals. Let $T, U$ be natural theories of consistency strength at least that of ``there are infinitely many Woodin cardinals''; then either $(\Pi^1_\omega)_T\subseteq(\Pi^1_\omega)_U$ or $(\Pi^1_\omega)_U\subseteq (\Pi^1_\omega)_T$.
\end{displayquote}
That is, at the level of $\Pi^1_\omega$ statements---i.e., statements about $\mathbb{R}$---all natural theories converge. Moreover, \textsc{steel's maxim} and its variants suggest that we ought to endorse the sentences on which they converge. Thus, the apparent pre-well-ordering of theories by consistency strength and other notions of proof-theoretic strength plays a central role in the search for new axioms in set theory.


\section{Ordinal Analysis}\label{section on ordinal analysis}

At this point it will be helpful to discuss two research programs that deal with questions of systematically reducing incompleteness. One is Gentzen's program of ordinal analysis and the other is Turing's program (significantly extended by Feferman) of completeness via iterated reflection principles. Both programs can be understood as reactions to the incompleteness theorems. The practitioners of ordinal analysis have attempted to reduce incompleteness by proving consistency statements in a systematic way, namely, from ever stronger transfinite induction principles. Turing's program, on the other hand, uses the second incompleteness theorem as the major engine for overcoming incompleteness; the practitioners of Turing's program try to systematically effect reductions in incompleteness by successively adding consistency statements to theories as new axioms. Both programs will play a major role in the approach to the ``great mystery'' discussed here.

Ordinal analysis was developed in the context of \emph{Hilbert's Program}, an early twentieth century research program pioneered by David Hilbert. To combat skepticism about the cogency of infinitary mathematics, Hilbert proposed to (i) axiomatize infinitary mathematics and (ii) prove the consistency of the axioms by finitary means. In 1931, Hilbert's program reached a major obstacle in the form of G\"{o}del's second incompleteness theorem \parencite{godel1931formal}. Indeed, it follows from G\"{o}del's theorem that if the principles of finitistic mathematics are codifiable in a reasonable axiomatic theory, then they do not prove their own consistency, much less the consistency of stronger theories. Thus, it is generally agreed that Hilbert's program failed. 

Gentzen was apparently undeterred, however. Not long after G\"{o}del proved the incompleteness theorems, Gentzen produced a consistency proof of arithmetic.

\begin{theorem}[\cite{gentzen1939}] If $\varepsilon_0$ is well-founded, then arithmetic is consistent.
\end{theorem}

The principles invoked in Gentzen's proof are nearly universally regarded as finitistically acceptable, except for the principle that the ordinal number $\varepsilon_0$ is well-founded. Gentzen's consistency proof marked the beginning of \emph{ordinal analysis}, a research program whereby similar consistency proofs have been discovered for a wide range of axiomatic theories. Developing such a consistency proof for a theory $T$ involves, among other things, determining the \emph{proof-theoretic ordinal} of $T$. Informally, the proof-theoretic ordinal of $T$ is the least ordinal $\alpha$ such that induction along $\alpha$ suffices to prove the consistency of $T$. Making this informal definition precise is no easy task, as we shall see.




The methods of ordinal analysis have been used to analyze many theories of interest, including subsystems of first-order arithmetic, second-order arithmetic, and set theory. Ordinal analysis has not yet reached the level of full second-order arithmetic. Nevertheless, it is reasonable to expect that the existing results are part of a general connection between ordinals and consistency. Do Gentzen-style methods suffice for proving the consistency of any axiomatic theory? The following result might \emph{seem} to suggest a positive answer ($\mathsf{TI}^\prec_{\Pi_1}$ is a sentence expressing the validity of induction for $\Pi_1$ predicates along $\prec$):
\begin{theorem}[\cite{kreisel1960number}] \label{scope of ordinal analysis} For any reasonable theory $T$ there is a primitive recursive presentation $\prec$ of an ordinal such that $\mathsf{PA}+\mathsf{TI}^\prec_{\Pi_0}\vdash \mathsf{Con}(T)$.\footnote{Even more dramatic versions of this result are developed in \cite{kreisel1976beweistheorie, kreisel1968survey}.}
\end{theorem}

This theorem is proved by showing that, given any true consistency statement $\varphi$, one can ``encode'' the truth of $\varphi$ into an ordinal notation system $\prec$. One would presumably not recognize the validity of $\Pi_1$ transfinite induction along such a notation system without knowing the truth of $\varphi$, so one could not use such a transfinite induction principle to prove $\varphi$. Thus, the epistemic value of this theorem is limited.

As we will see, this is one version of a pervasive problem known as the \emph{canonicity problem}. The ordinal notations that are devised to prove Theorem \ref{scope of ordinal analysis} are not \emph{natural} notations. They form a notation system that one would introduce only in an ad-hoc manner to solve a problem in proof theory. Does the distinction between the pathological notation systems and the natural notation systems reflect some intrinsic mathematical properties of the notation systems? If one had a convincing definition of ``natural,'' one might conjecture that, for any reasonable theory $T$, there is a \emph{natural} presentation $\prec$ of a recursive ordinal such that $\mathsf{PA}+\mathsf{TI}^\prec_{\Pi_1}\vdash \mathsf{Con}(T)$. However, at present there is no convincing evidence that it is possible to precisely define the ``natural'' notation systems.

The canonicity problem also makes it difficult to define ``proof-theoretic ordinal,'' as suggested earlier. One might try to define the proof-theoretic ordinal of $T$ as the least ordinal $\alpha$ such that induction along $\alpha$ (along with finitary methods) suffices to prove the consistency of $T$. The problem with this definition is that, in formalized theories, transfinite induction principles are stated in terms of \emph{presentations} of ordinals, not the ordinals themselves. \cite{kreisel1976beweistheorie} shows that it is always possible to prove the consistency of a reasonable theory by induction along a sufficiently pathological presentation of $\omega$. Conversely, \cite{beklemishev2001another} shows that, for any reasonable theory $T$ and recursive ordinal $\alpha$, there is a sufficiently pathological presentation of $\alpha$ such that transfinite induction along that presentation (along with finitary methods) does not prove the consistency of $T$. When one restricts one's attention to ``natural'' presentations of ordinals, the na\"ive definition seems to work, but the current state of affairs is vexing and unsatisfactory.

At present, many conflicting definitions of the ``proof-theoretic ordinal'' of a theory have been proposed; these definitions often coincide in crucial cases. Perhaps the most common definition is this: The proof-theoretic ordinal of a theory $T$ is the supremum of the order-types of the primitive recursive well-orderings whose well-foundedness is provable in $T$. This is sometimes called the $\Pi^1_1$ ordinal of a theory. Of course, this notion is useful only for measuring the proof-theoretic strength of theories in which well-foundedness is expressible, so it does not apply, e.g., to fragments of first-order arithmetic. Moreover, it is a somewhat coarse notion of strength, since $\Pi^1_1$ ordinals are invariant under the addition of true $\Sigma^1_1$ sentences to the object theory. Nevertheless, the $\Pi^1_1$ ordinal of a theory is an \emph{ordinal}, not a presentation of an ordinal, so this is a somewhat robust notion.

\section{Turing--Feferman Progressions}\label{turing-progressions}

\subsection{Turing's Completeness Theorem}
The second incompleteness theorem suggests a method for dealing with the first incompleteness theorem. According to Turing, if we endorse the axioms of a reasonable axiomatic theory $T$, then there is a principled way of extending $T$, namely, by adopting $T$'s consistency statement as an axiom. Of course, if one adopts the statement $\mathsf{Con}(T)$ as an axiom, the statement $\mathsf{Con}(T+\mathsf{Con}(T))$ remains unprovable. However, there is a principled way of resolving this problem, namely, by adopting $\mathsf{Con}(T+\mathsf{Con}(T))$ as an axiom. Of course, this engenders a new problem, but it just as easily engenders a new solution, whence this process can be iterated \emph{ad infinitum}. Given presentations of recursive ordinals in the language of arithmetic, one can even extend this process into the effective transfinite.

Na\"ively, we might try to define the iterations of consistency over a theory $T$ as follows:
\begin{itemize}
\item $T_0:=T$
\item $T_{\alpha+1} := T_\alpha + \mathsf{Con}(T_\alpha)$
\item $T_\lambda :=\bigcup_{\alpha<\lambda}T_\alpha$ for $\lambda$ a limit.
\end{itemize}
However, such a definition does not even pin down the theory $T_{\omega+1}$. According to the definition, $T_{\omega+1}$ is just $T_\omega+\mathsf{Con}(T_\omega)$, but for $\mathsf{Con}(T_\omega)$ to be a statement of arithmetic we must have some effective presentation of $T_\omega$, and there are many choices for the latter. Accordingly, when one defines iterated consistency statements, one must first fix an ordinal notation system $\prec$. One can then define the iterations of consistency along $\prec$ within arithmetic via G\"{o}del's fixed point lemma. Let me briefly indicate how this is done; thenceforth I will only give the intuitive ``definitions'' by transfinite recursion, trusting the reader to fill in the details.\footnote{For a discussion of this technique and comparisons with other techniques see the proof of Theorem 3.8 in \cite{feferman1962transfinite}.}


We fix a formula $Ax(x,T)$ which formalizes that $x$ is an axiom of $T$. We now want to define a formula $Ax(x,\alpha,T)$ that defines the axioms of the theory $T_\alpha$. We define the formula as the following fixed point:
$$\mathsf{PA} \vdash Ax(x,\alpha,T) \leftrightarrow \big(Ax(x,T) \vee \exists \beta\prec \alpha \; x = \ulcorner\mathsf{Con}(T_\beta)\urcorner \big).$$
Put a bit more explicitly:
{\footnotesize
$$\mathsf{PA} \vdash Ax(x,\alpha,T) \leftrightarrow \Big(Ax(x,T) \vee \exists \beta\prec \alpha \; x=\ulcorner \forall z \mathsf{Con}\big(\bigwedge \{ w \mid w\leq z \wedge  Ax(w,\beta,T) \}\big)\urcorner\Big).$$
}
Note that the three-place formula $Ax(\cdot,\cdot,\cdot)$ appears on the right but only in G\"odel quotes, whence this is a legitimate application of the fixed point lemma.

Does iteratively endorsing consistency statements in this manner effect a significant reduction in incompleteness? The following result may \emph{seem} to suggest a positive answer:

\begin{theorem}[\cite{turing1939systems}] For any true $\Pi_1$ sentence $\varphi$, there is a presentation $\dot\alpha$ of $\omega+1$ such that $\mathsf{PA}_{\dot\alpha}\vdash\varphi$.
\end{theorem}

This result is known as \emph{Turing's Completeness Theorem}. At first glance this theorem seems epistemically significant: We can come to know the truth of any $\Pi_1$ statement $\varphi$ simply by iterating consistency statements. How are the consistency statements used in the proof of $\varphi$? The disappointing response is that they are not used at all. Instead, the truth of $\varphi$ is encoded into a non-standard description of the base theory $\mathsf{PA}$, which becomes available at iteration $\omega+1$. Given this non-standard description of $\mathsf{PA}$, discerning what the theory $\mathsf{PA}_{\dot\alpha}$ is committed to \emph{requires} knowing the truth-value of $\varphi$. This is another appearance of the canonicity problem. As far as I know, this was the initial appearance of the problem, and Turing was the first to identify it.

Turing was careful not to overstate the epistemic significance of his theorem. Indeed, he wrote:
\begin{displayquote}[\cite{turing1939systems}, p. 213]
This completeness theorem as usual is of no value. Although it shows, for instance, that it is possible to prove Fermat’s last theorem\dots (if it is true) yet the truth of the theorem would really be assumed by taking a
certain formula as an ordinal formula.
\end{displayquote}

Nevertheless, one might hope that iterating $\mathsf{Con}$ along \emph{natural} ordinal notations would also suffice to prove any true $\Pi_1$ statement. Once again, if we had a precise definition of ``natural'' ordinal notation systems, we would be able to formulate a precise conjecture; however, no such characterization of ``natural'' ordinal notation systems is currently available.

Turing discussed the possibility of using iterated consistency progressions to measure the strength of $\Pi_1$ statements. He wrote:
\begin{displayquote}[\cite{turing1939systems} p. 200]
We might also expect to obtain an interesting classification of number-theoretic theorems according to ``depth''. A theorem which required an ordinal $\alpha$ to prove it would be deeper than one which could be proved by the use of an ordinal $\beta$ less than $\alpha$. However, this presupposes more than is justified.
\end{displayquote}
Na\"ively, one approach would be to assign to each true $\Pi_1$ sentence $\varphi$ the least ordinal $\alpha$ such that $\mathsf{PA}_\alpha\vdash\varphi$ as its measure of complexity. One could then extend this measure of complexity to theories: the complexity of $T$ is the least ordinal such that $\mathsf{PA}_\alpha\vdash\varphi$ for each $\Pi_1$ theorem $\varphi$ of $T$. Clearly, Turing's completeness theorem shows that iterations of consistency depend on presentations of ordinals and not just on ordinals themselves. The possibility remains, however, that one could prove informative results of this sort for large swathes of interesting theories by antecedently fixing some natural notation system. Beklemishev has recently pursued this possibility, as we will discuss in \textsection \ref{Ordinal analysis via iterated reflection}.

\subsection{Reflection Principles}\label{reflection-definitions}
Turing's work was later pursued and greatly extended by Feferman \parencite{feferman1962transfinite}. Feferman shifted the focus from iterated consistency to iterated reflection principles of other sorts. We pause here to describe these reflection principles.

Feferman's focus was on the full uniform reflection schema, denoted $\mathsf{RFN}(T)$. To discuss this schema, we must review some notation. In what follows, $\overline{m}$ is the standard numeral of the number $m$ and $\ulcorner\varphi\urcorner$ is the G\"odel number of the formula $\varphi$. The expression $\ulcorner \varphi(\dot{x_1},\dots,\dot{x_n})\urcorner$ denotes an elementary definable term for the function taking $k_1,\dots,k_n$ to the G\"odel number $\ulcorner\varphi(\overline{k_1},\dots, \overline{k_n} )\urcorner$ of the result of substituting the numerals $\overline{k_1},\dots, \overline{k_n}$ for the variables $x_1,\dots,x_n$ in $\varphi$.

The full uniform reflection principle for $T$ is the schema:
$$\mathsf{RFN}(T):= \{ \forall x_1\dots\forall x_n \Big(\mathsf{Pr}_T\big(\ulcorner\varphi(\dot{x_1},\dots,\dot{x_n})\urcorner\big) \to \varphi (x_1,\dots,x_n) \Big) \mid \varphi \in\mathcal{L}_A\}.$$


We sometimes consider restrictions of the full uniform reflection schema to various complexity classes. One such schema is the restriction of the uniform reflection principle to $\Sigma_n$ formulas, i.e:
$$\{ \forall x_1\dots\forall x_n \Big(\mathsf{Pr}_T\big(\ulcorner\varphi(\dot{x_1},\dots,\dot{x_n})\urcorner\big) \to \varphi (x_1,\dots,x_n) \Big) \mid \varphi \in \Sigma_n\}.$$
In fact, this principle is axiomatizable over weak systems\footnote{$\mathsf{EA}$ suffices. For a definition of $\mathsf{EA}$, see \textsection \ref{Ordinal analysis via iterated reflection}.} by a single sentence, namely:
$$\mathsf{RFN}_{\Sigma_n}(T):= \forall \varphi \in \Sigma_n \big( \mathsf{Pr}_T(\ulcorner\varphi\urcorner)\rightarrow \mathsf{True}_{\Sigma_n}(\ulcorner\varphi\urcorner) \big).$$
Likewise, the schema of uniform $\Pi_n$ reflection:
$$\{ \forall x_1\dots\forall x_n \Big(\mathsf{Pr}_T\big(\ulcorner\varphi(\dot{x_1},\dots,\dot{x_n})\urcorner\big) \to \varphi (x_1,\dots,x_n) \Big) \mid \varphi \in \Pi_n\}.$$
also has a finite axiomatization:
$$\mathsf{RFN}_{\Pi_n}(T):= \forall \varphi \in \Pi_{n} \big( \mathsf{Pr}_T(\ulcorner\varphi\urcorner)\rightarrow \mathsf{True}_{\Pi_{n}}(\ulcorner\varphi\urcorner) \big).$$

The last reflection principle that we introduce is that of $n$-consistency. We say that a theory is $n$-consistent if it is consistency with the true $\Pi_n$ theory of arithmetic. We can express that a theory $T$ is $n$-consistent as follows:\footnote{It is conventional to omit corner quotes when designating theories like $T+\varphi$ inside $\mathsf{Con}$ even though the formula literally mentions the G\"odel number of $\varphi$. See, e.g., \cite{beklemishev2004algebras}.}
$$ n\mathsf{Con}(T):= \forall \varphi \in \Pi_{n}\big( \mathsf{True}_{\Pi_{n}}(\ulcorner\varphi\urcorner) \to \mathsf{Con}(T+\varphi )\big).$$

In fact, over $\mathsf{EA}$, the statements $\mathsf{RFN}_{\Sigma_n}(T)$, $\mathsf{RFN}_{\Pi_{n+1}}(T)$, and $ (n+1)\mathsf{Con}(T)$ are all provably equivalent.

\subsection{Feferman's Completeness Theorem}
Feferman's theorem concerns the iteration of the full uniform reflection schema. After fixing a presentation $\prec$ of a recursive ordinal, we can then define iterations of reflection over a theory $T$ so as to satisfy the following conditions:
\begin{itemize}
\item $T^{\mathsf{RFN}}_0:=T$
\item $T^{\mathsf{RFN}}_{\alpha+1} := T^{\mathsf{RFN}}_\alpha + \mathsf{RFN}(T^{\mathsf{RFN}}_\alpha)$
\item $T^{\mathsf{RFN}}_\lambda :=\bigcup_{\alpha<\lambda}T^{\mathsf{RFN}}_\alpha$ for $\lambda$ a limit.
\end{itemize}

Feferman proved that, merely by iterating the uniform reflection schema along presentations of recursive ordinals, one can accrue resources sufficient, not only for proving any true $\Pi_1$ statement, but for proving \emph{any} true arithmetical statement.

\begin{theorem}[\cite{feferman1962transfinite}] For any true arithmetical sentence $\varphi$, there is a presentation $\dot\alpha$ of an ordinal $\alpha<\omega^{\omega^\omega+1}$ such that $\mathsf{PA}^{\mathsf{RFN}}_{\dot\alpha}\vdash\varphi$.
\end{theorem}

This result is known as \emph{Feferman's Completeness Theorem}. In the proof of this theorem, Feferman used Turing's technique of encoding the truth of statements into presentations of recursive ordinals, among other things. Feferman also proved that there are paths $P$ through Kleene's $\mathcal{O}$ such that iterating uniform reflection along $P$ is arithmetically complete \parencite{feferman1962transfinite}. However, no such path is even $\Sigma^1_1$ definable \parencite{feferman1962incompleteness}. That is, identifying such a path is \emph{more difficult} than identifying the set of arithmetical truths. Thus, the epistemic significance of these results---as Feferman emphasized---is limited. 

\section{Ordinal Analysis via Iterated Reflection}\label{Ordinal analysis via iterated reflection}

In recent years there has been interest in the interface between ordinal analysis and Turing progressions. This research has been motivated by a number of the drawbacks of the research discussed in the previous two sections.

As discussed in \textsection \ref{section on ordinal analysis}, the standard notion of the $\Pi^1_1$ ordinal of a theory is (i) only applicable to theories in which ``well-foundedness'' is expressible and (ii) insensitive to the true $\Sigma^1_1$ consequences of a theory. Building on \cite{schmerl1979fine}, \cite{beklemishev2003proof} introduced the notion of the $\Pi_1$ ordinal of a theory. This notion of ``ordinal analysis'' is both (i) suitable for subsystems of first-order arithmetic and (ii) sensitive to the $\Pi_1$ consequences of theories. Iterations of consistency in the style of Turing play a central role in the definition of $\Pi_1$ ordinals.

First we fix a base theory: Beklemishev uses $\mathsf{EA}$, a weak subsystem of arithmetic. $\mathsf{EA}$ proves the totality of exponentiation but not superexponentiation; $\mathsf{EA}$ is just strong enough to carry out arithmetization of syntax in the standard way. We then fix some natural ordinal notation system and an elementarily defined extension $T$ of $\mathsf{EA}$. We define the iterations of $T$ so as to satisfy the following conditions:
\begin{enumerate}
\item $\mathsf{EA}_0:=\mathsf{EA}$
\item $\mathsf{EA}_{\alpha+1}:=\mathsf{EA}_\alpha+\mathsf{Con}(\mathsf{EA}_\alpha)$
\item $\mathsf{EA}_\lambda : =\bigcup_{\alpha<\lambda}\mathsf{EA}_\alpha$ for $\lambda$ a limit.
\end{enumerate}
Given a target theory $T$, the $\Pi_1$ ordinal of $T$ (relative to the base theory $\mathsf{EA}$ and notation system $\prec$) is defined as follows:
$$|T|_{\Pi_1}:=\mathsf{sup}\{\alpha:\mathsf{EA}_\alpha \subseteq T\}$$
This definition yields interesting information for theories that are conservatively approximated by iterations of consistency over $\mathsf{EA}$; that is, for any $T$ such that $T\equiv_{\Pi_1} \mathsf{EA}_\alpha$ where $\alpha=|T|_{\Pi_1}$. Knowing that $T$ is so approximated is useful, because there are elegant equations, first discovered in \cite{schmerl1979fine} but refined in \cite{beklemishev2003proof}, which spell out conservation relations between iterated reflection principles. $\Pi_1$ ordinals are well-defined for theories of first-order arithmetic in which well-foundedness is not directly expressible. The definition is also sensitive to $\Pi_1$ sentences. For instance, whereas the $\Pi_1$ ordinal of $\mathsf{PA}$ is $\varepsilon_0$, the $\Pi_1$ ordinal of $\mathsf{PA}+\mathsf{Con}(\mathsf{PA})$ is $\varepsilon_0\times 2$.

Note that the definition of $\Pi_1$ ordinals achieves Turing's goal of providing a metric of strength for theories in terms of the number of iterations of consistency required to capture their $\Pi_1$ theorems. In the cases that have been investigated, it turns out that Turing's metric of strength is the same one provided by ordinal analysis. That is, the definition of $\Pi_1$ ordinals coincides with other standard definitions of proof-theoretic ordinal in standard cases. To provide such a metric of strength, we must fix (in advance) a natural ordinal notation system. Moreover, the definition works only for theories whose $\Pi_1$ fragments can be conservatively approximated by iterated consistency statements. This is apparently a feature of ``natural theories'' but not of all theories. This approach to ordinal analysis was a major inspiration for the approach to the main topic of this paper discussed in \textsection \ref{operators}--\ref{reflection-and-ordinals}.

These ideas are closely related to Beklemishev's approach to the canonicity problem for ordinal notation systems \parencite{beklemishev2004provability, beklemishev2005reflection}. I will briefly describe the terms that make up Beklemishev's ordinal notation system for analyzing $\mathsf{PA}$ and its fragments. The terms are generated by the constant symbol $\top$ and the function symbols $\Diamond_0(\cdot), \Diamond_1(\cdot), \Diamond_2(\cdot),$ etc. The ordering on these terms is defined in terms of a variable-free polymodal logic $\mathsf{GLP}_0$.\footnote{Alternatively, we may use the computationally simpler Reflection Calculus $\mathsf{RC}$.} In particular:
$$\alpha \prec \beta \Leftrightarrow \mathsf{GLP}_0 \vdash \beta \rightarrow \Diamond_0 \alpha.$$


Let me emphasize one aspect of this term system that is relevant to present concerns. First, we define a translation $\star$ from the variable-free polymodal language into the language of arithmetic such that:
\begin{itemize}
    \item $\top^\star := 0=0$
    \item $\Diamond_n(\varphi)^\star := n\mathsf{Con}(T+\varphi^\star)$
    \item $\star$ commutes with Boolean connectives.
\end{itemize}
For every $\varphi$ of the variable-free polymodal language, $\mathsf{GLP}_0$ proves $\varphi$ if and only if $\mathsf{EA}$ proves $\varphi^\star$.\footnote{We state things here in terms of a single translation $\star$ since we are interested only in the variable-free system $\mathsf{GLP}_0$. If we introduce propositional variables and look at all translations in the usual way, $\mathsf{GLP}$ is still sound with respect to its arithmetical interpretation but completeness is an open problem.} So the ordering on terms corresponds with consistency strength over $\mathsf{EA}$ in the following sense:
$$ \text{For terms $\alpha$ and $\beta$, $\alpha \prec \beta$ if and only if $\mathsf{EA}\vdash \alpha^\star \to \Diamond_0 \beta^\star$.}$$
This notation system is well-suited to ordinal analysis via iterated reflection. $\mathsf{PA}$ can be conservatively approximated over $\mathsf{EA}$ by the interpretations of the terms of the notation system.

The connection with the canonicity problem is this: If these reflection principles are regarded as a \emph{canonical} means of specifying the theory $\mathsf{PA}$, then the ordinal notation system and not just the ordinal has been extracted from a \emph{canonical} presentation of $\mathsf{PA}$. The results discussed in this paper support the idea that reflection principles \emph{are} canonical. Moreover, these results support the idea that natural theories can be approximated by iterating reflection principles. It remains to be seen whether the canonicity of reflection principles can be established in a way that bears on the problem of canonical ordinal notations. That is, the results discussed in this paper, though they show that reflection principles are in \emph{some} sense canonical, seem too coarse grained to isolate Beklemishev's notation system as canonical.\footnote{The source of coarseness is that we will show that the standard reflection principles are canonical only \emph{up to} a certain equivalence relation.}

\section{Turing Degree Theory}\label{turing-degrees}

The pre-well-ordering of natural theories is paralleled by phenomena in Turing degree theory. Let's start with a discussion of \emph{Post's Problem}. Post asked whether there is a recursively enumerable degree $d$ such that $0 <_T d <_T 0'$. This question was answered positively by Friedberg and Muchnik. However, they constructed such intermediate degrees by deploying ad hoc methods; the degrees that they construction are widely regarded as ``unnatural.'' Is there a \emph{natural} recursively enumerable degree $d$ such that $0 <_T d <_T 0'$? Once again, it is not entirely clear how to address this problem. It is not a wholly mathematical problem since there is no precise mathematical definition of ``natural'' Turing degrees.

Sacks clarified this situation by asking whether there is a degree-invariant solution to Post's Problem. Let me briefly explain.
\begin{definition}
    A function $W : 2^\omega \to 2^\omega$ is a \emph{recursively enumerable operator} if there is an $e\in\mathbb{N}$ such that, for each $A$, $W(A)=W_e^A$.
\end{definition}
\begin{definition}
    A function $f:2^\omega \to 2^\omega$ is \emph{degree-invariant} if for all reals $A$ and $B$, $A\equiv_T B$ implies $f(A)\equiv_T f(B)$.
\end{definition}
One oft-noted feature of natural Turing degrees is that their definitions \emph{relativize}. Relativizing the definition of a Turing degree yields a degree-invariant function on the reals. For instance, relativizing the definition of $0'$ yields the Turing jump $\lambda X.X'$. We are now in a position to state Sacks' question.
\begin{question}[\cite{sacks1966degrees}]
    Is there a degree-invariant solution to Post's Problem? That is, is there a degree-invariant recursively enumerable operator $W$ such that for every real $A$, $A<_T W(A) <_T A'$?
\end{question}
This question is still unresolved. A negative answer would explain why the only natural recursively enumerable degrees are $0$ and $0'$. Relativizing the definition of any natural recursively enumerable degree should also yield a degree-invariant recursively enumerable operator. Thus, if there are no degree-invariant recursively enumerable operators strictly between the identity and the jump, it suggests that there are no natural recursively enumerable degrees other than $0$ and $0'$.

Sacks' question is actually a special case of a more wide-reaching question that has motivated a considerable amount of research into the Turing degrees. The Turing degrees are not linearly ordered by $<_T$. That is, there are distinct degrees $a$ and $b$ such that $a\nleq_Tb$ and $b\nleq_Ta$. The Turing degrees are also ill-founded. That is, there are infinite sequences $(a_n)_{n<\omega}$ such that $a_k>_Ta_{k+1}$ for all $k$. These two results mean that it is neither possible to compare nor to rank Turing degrees in general. This behavior is exhibited below $0'$, but it is in fact realized throughout the Turing degrees.

The degrees that exhibit these pathological properties have been constructed using ad-hoc recursion-theoretic techniques, like the priority method. When one restricts one's attention to \emph{natural} Turing degrees, the resulting structure is a well-ordering.
$$ 0, 0', 0'', \dots, 0^\omega, \dots, \mathcal{O}, \dots, 0^\sharp, \dots $$
This state of affairs is remarkably similar to the state of affairs on the proof-theoretic side. Once again, it is not entirely clear how to address this problem. It is not a wholly mathematical problem since there is no precise mathematical definition of ``natural'' Turing degrees.

As already mentioned, relativizing the definition of natural Turing degrees yields a degree-invariant function on the reals. Just as relativizing the definition of $0'$ yields the Turing jump $\lambda X.X'$, relativizing the definition of Kleene's $\mathcal{O}$ yields the hyperjump $\lambda X.\mathcal{O}^X$, relativizing the definition of $0^\sharp$ yields the sharp function $\lambda X.X^\sharp$, and so on. The identity function and constant functions are also degree-invariant. Are these all the degree-invariant functions that there are? There are indeed many others. Indeed, for any real $A$, we can concoct a degree-invariant function $\lambda X. A \oplus X'$. Yet this function is clearly just a modification of the Turing jump, and if we restrict our attention to the Turing degrees above $A$, this function behaves \emph{exactly} as the Turing jump does. Other examples of degree-invariant functions also seem to be mere variations of the functions already described. So one might suspect that---modulo a suitable equivalence relation---the identity, the constant functions, and generalizations of the Turing jump \emph{are} the only degree-invariant functions.

Martin proposed a classification of the degree-invariant functions that makes sense of this informal suspicion. Martin's Conjecture classifies the degree-invariant functions in terms of their behavior \emph{almost everywhere} in the sense of a measure known as \emph{Martin Measure}.\footnote{For an accessible overview of the conjecture, see \cite{montalban2019martin}.} To state the definition of Martin Measure, recall that a \emph{cone} in the Turing degrees is any set of the form $\{a: a\geq_T b\}$. Assuming $\mathsf{AD}$, Martin proved that every degree-invariant set of reals either contains a cone or is disjoint from a cone. Moreover the intersection of countably many cones contains a cone. Thus, assuming $\mathsf{AD}$, the function
\[ \mu(A)= \begin{cases} 
      1 & \textrm{ if $A$ contains a cone} \\
      0 & \textrm{ if $A$ is disjoint from a cone} 
   \end{cases}
\]
is a countably additive measure on the $\sigma$ algebra of degree-invariant sets. This measure is called \emph{Martin Measure}. In the statement of Martin's Conjecture, \emph{almost everywhere} means \emph{almost everywhere with respect to Martin Measure}. If $f(x)\leq_T g(x)$ almost everywhere (in the sense of Martin Measure), then we say that $f$ is \emph{Martin below} $g$ and write $f\leq_M g$.

\begin{conjecture}[Martin]
Assume $\mathsf{ZF}+\mathsf{DC}+\mathsf{AD}$. Then
\begin{enumerate}
\item[I.] If $f:2^\omega\rightarrow2^\omega$ is degree-invariant, and $f$ is not increasing almost everywhere then $f$ is constant almost everywhere.
\item[II.] $\leq_M$ pre-well-orders the set of degree-invariant functions that are increasing almost everywhere. If $f$ has $\leq_M$ rank $\alpha$, then $f'$ has $\leq_M$ rank $\alpha+1$, where $f'(x)=f(x)'$ for all $x$.
\end{enumerate}
\end{conjecture}

Martin's Conjecture is stated under the Axiom of Determinacy, $\mathsf{AD}$. Though one can refute $\mathsf{AD}$ in $\mathsf{ZFC}$, one cannot refute restrictions of $\mathsf{AD}$ to sets of reals with sufficiently simple definitions.\footnote{That is, assuming that $\mathsf{ZFC}$ is consistent.} Indeed, using $\mathsf{ZFC}$, one can prove Borel Determinacy (the restriction of the Axiom of Determinacy to Borel sets). Assuming appropriate large cardinal axioms, one can even prove that $L(\mathbb{R})$ satisfies $\mathsf{AD}$. That is, as we make stronger and stronger assumptions, we can prove restrictions of determinacy to the ``definable'' sets for ever more liberal notions of definability. Hence, Martin's Conjecture is often understood as a conjecture about the ``definable'' functions. Understood this way, Martin's Conjecture roughly states that the only definable degree-invariant functions (up to almost everywhere equivalence) are constant functions, the identify function, and iterates of the Turing Jump.

Though Martin's Conjecture is presently open, many informative partial results and special cases are known. Let me briefly mention some.\footnote{Also worth mentioning, though it is less relevant to present concerns, is Martin's Conjecture for uniformly degree-invariant functions. Part I is proved in \cite{slaman1988definable} and Part II is proved in \cite{steel1982classification}.}

\begin{definition}
    A function $f$ is \emph{order-preserving} if $A\leq_T B$ implies $f(A)\leq_T f(B)$. 
\end{definition}




\begin{theorem}[\cite{lutz2023part}]
Part I of Martin's Conjecture holds for order-preserving functions.
\end{theorem}

\begin{theorem}[\cite{slaman1988definable}]
Part II of Martin's Conjecture holds for order-preserving Borel functions.
\end{theorem}

There are other results that speak to this phenomenon of well-foundedness. For instance, Steel proved that there cannot be any ``simple'' descending sequences in the Turing jump hierarchy, in the following sense:

\begin{theorem}[\cite{steel1975descending}]\label{steel-descending}
Let $P\subset \mathbb{R}^2$ be arithmetic. Then there is no sequence $(x_n)_{n<\omega}$ such that for every $n$,
\begin{itemize}
\item[(i)] $x_n\geq_T x'_{n+1}$ and
\item[(ii)] $x_{n+1}$ is the unique $y$ such that $P(x_n,y)$.
\end{itemize}
\end{theorem}

We will explore the analogy between the recursion-theoretic and proof-theoretic well-ordering phenomena throughout the rest of this paper. In particular, we will look at proof-theoretic results inspired by Martin's Conjecture and its partial realizations.


\section{Interlude}\label{interlude}

My goal is to bring precision to the question of the pre-well-ordering of natural theories and to offer (at least partial) solutions. Before summarizing this research, I would like to mention three themes that will be interwoven.

The first theme is that many natural theories can be axiomatized by reflection principles over natural base theories. Indeed, the fragments of natural theories corresponding to different syntactic complexity classes can often be conservatively approximated by iterated reflection principles of the appropriate complexity class. 

The second theme is that ordinal analysis is an expression of the well-ordering phenomenon. For different notions of ``proof-theoretic strength,'' there are corresponding notions of ``proof-theoretic ordinals.''  Insofar as the proof-theoretic ordinal afforded by some definition of ``proof-theoretic ordinal'' measures the proof-theoretic strength of theories, the attendant method of ordinal analysis well-orders the theories within its ken according to that notion of strength. The heuristic that we will try to vindicate in this paper is the following: To calculate the proof-theoretic ordinal of a theory $T$ is to determine $T$'s rank in the hierarchy of natural theories ordered by proof-theoretic strength.

The third theme is that reflection principles play a role in proof theory analogous to the role jumps play in recursion theory. Just as natural Turing degrees are apparently equivalent to ordinal iterates of the Turing jump, natural theories are apparently equivalent to ordinal iterates (along natural presentations of well-orderings) of reflection principles. Just as natural Turing degrees can be obliquely studied in terms of jumps, natural axiomatic theories can be obliquely studied in terms of reflection principles.


\section{Analogues of Martin's Conjecture}\label{operators}

Ordinal analysis is a means of well-ordering axiom systems, so one might expect some connections between ordinal analysis and the well-ordering phenomenon. Recall that Beklemishev's method (see \textsection \ref{Ordinal analysis via iterated reflection}) for calculating the proof-theoretic ordinal of a theory $T$ involves demonstrating that $T$ can be approximated over a weak base theory by a class of formulas that are well-behaved. In particular, the $\Pi_1$ fragments of natural theories are often proof-theoretically equivalent to iterated consistency statements over a weak base theory, whence these theories are amenable to ordinal analysis. Why are the $\Pi_1$ fragments of natural theories proof-theoretically equivalent to iterated consistency statements?

Our approach to this problem is inspired by Martin's Conjecture and similar problems in recursion theory. In this section we will formulate and discuss proof-theoretic versions of these problems. Let's fix some notation. Throughout this section we will let $T$ be a sound (i.e., true in the standard model), recursively axiomatized same-language extension of elementary arithmetic. We use the following notation frequently:
\begin{definition}
    $[\varphi]_T$ is the equivalence class of $\varphi$ modulo $T$ provable equivalence, i.e., $[\varphi]_T = \{\psi: T\vdash \varphi \leftrightarrow \psi\}$.
\end{definition}
Finally, in this section, whenever we write ``true'' we mean ``true according to the standard interpretation $\mathbb{N}$ of first-order arithmetic.''

Now, let's consider an analogue of Post's Problem. Are there theories whose deductive strength lies strictly between that of $T$ and $T+\mathsf{Con}(T)$? The answer is positive. Indeed, there is a true $T$-independent $\Pi_1$ sentence $\varphi$ such that $T+\mathsf{Con}(T)\vdash \varphi$ but $T+\varphi\nvdash \mathsf{Con}(T)$. Indeed, by the second incompleteness theorem $T+\neg \mathsf{Con}(T)$ is consistent. Thus, by Rosser's Theorem there is some sentence that is neither provable nor refutable in $T+\neg \mathsf{Con}(T)$. In particular, Rosser's theorem shows that the sentence $R$ is $T+\neg \mathsf{Con}(T)$-independent, where:
$$T\vdash  R \leftrightarrow  \forall x \big(\mathsf{Prf}_{T+\neg \mathsf{Con}(T)}(x,\ulcorner R \urcorner) \rightarrow \exists y<x \; \mathsf{Prf}_{T+\neg \mathsf{Con}(T)}(y,\ulcorner\neg R\urcorner) \big) .$$
$R$ ``says'' that all $T+\neg \mathsf{Con}(T)$ proofs of $R$ (if there are any) are preceded by $T+\neg \mathsf{Con}(T)$ proofs of $\neg R$.
An elementary argument then shows that $T+ \big(\mathsf{Con}(T) \vee R \big) $ has deductive strength strictly between  $T$ and $T+\mathsf{Con}(T)$.

The theory $T+ \big(\mathsf{Con}(T) \vee R \big)$ has the properties we were looking for, but it is a highly unnatural theory. Nobody would have proposed $\mathsf{Con}(T) \vee R$ as a principle of \emph{arithmetic}; it was defined in an unusual way to exhibit quirky \emph{meta-mathematical} features. Is there any way of explaining the unnaturalness of this formula in a mathematically exact way?


Recall that $R$ was defined in terms of the theory $T+\neg\mathsf{Con}(T)$. Indeed, Rosser proved that for \emph{any} reasonable theory $U$ there is a $U$-independent sentence $R_U$ where:
$$U\vdash  R_U \leftrightarrow  \forall x \big(\mathsf{Prf}_{U}(x,\ulcorner R_U \urcorner) \rightarrow \exists y<x \; \mathsf{Prf}_{U}(y,\ulcorner\neg R_U\urcorner) \big) .$$
Note that $R_U$ depends on some seemingly arbitrary details about the way that the theory $U$ is presented. In particular, $R$ makes a claim about the order in which proofs from the theory $U$ are encoded by natural numbers. Yet there are \emph{many} different ways of numerating proofs from the theory $U$. Some of these differences are merely differences in how $U$ itself is presented or described (for instance, the order in which the axioms of $U$ are listed). So, for instance, there could be two sentences $\varphi$ and $\psi$ such that $T\vdash \varphi \leftrightarrow \psi$ but $T\nvdash R_{T+\varphi} \leftrightarrow R_{T+\psi}.$

All of this is to say that the function $\varphi \mapsto R_{T+\varphi}$ is not extensional over $T$.
\begin{definition}
A function $f:
\mathbb{N}\to \mathbb{N}$ is \emph{extensional} over a theory $T$ if for all $\varphi$ and $\psi$:
$$T \vdash \varphi \leftrightarrow \psi \Longrightarrow T \vdash f(\varphi) \leftrightarrow f(\psi).$$
\end{definition}
\begin{remark}
    Note that the previous definition technically involves an abuse of notation. Since $f$ is a function from $\mathbb{N}$ to $\mathbb{N}$, it should map the \emph{G\"odel numbers} of $\varphi$ and $\psi$ to the \emph{G\"odel numbers} of other formulas. If $\theta_0,\theta_1,\dots$ is a fixed G\"odel numbering of the formulas of arithmetic, then the following gives an exact description of the condition of extensionality: 
$$T\vdash \varphi \leftrightarrow \psi \Longrightarrow T\vdash \theta_{f (\ulcorner\varphi\urcorner)} \leftrightarrow \theta_{f (\ulcorner\psi\urcorner)}.$$
We find that this particular abuse of notation is a particularly aid for intuition and we henceforth use it without comment.
\end{remark}
The output of an extensional function with input $\varphi$ depends only on $[\varphi]_T$. The identity function and constant functions are all extensional. So is the consistency operator $\varphi \mapsto \mathsf{Con}(T+\varphi)$ and its iterates into the effective transfinite. We introduce the following notation to reduce clutter in what follows:
$$\mathsf{Con}_T(\ulcorner \varphi \urcorner) := \mathsf{Con}(T+\varphi). $$

We will mostly be concerned with functions that are not only extensional but also monotone over $T$.
\begin{definition}
A function $f:
\mathbb{N}\to \mathbb{N}$ is \emph{monotone} over a theory $T$ if for all $\varphi$ and $\psi$:
$$T \vdash \varphi \to \psi \Longrightarrow T \vdash f(\varphi) \to f(\psi).$$
\end{definition}

Note that all monotone functions are extensional but not vice-versa. The goal of this approach is to show that the consistency operator and its iterates are canonical monotone functions.

\subsection{A Positive Theorem}

Let's begin with a positive result.\footnote{These results are not presented in the order that they were proved or published. I have rearranged them to increase readability. The dates in the statements of the theorems reflect the order in which they were published.} This result shows that sufficiently constrained recursive monotone functions must coincide with the consistency operator \emph{in the limit} within the ultrafilter of true sentences. Note that ``in the limit'' plays a role analogous to that which ``a.e.\ in the sense of Martin measure'' does in Martin's Conjecture. Here is how we precisely define ``in the limit'':
\begin{definition}
    A \emph{cone} is any set of sentences of the form $\{\varphi \mid T+\varphi \vdash \psi\}$, where $\psi$ is a sentence. A \emph{true cone} is a cone that contains a sentence that is true according to the standard interpretation of the language of arithmetic.
\end{definition}

One may think of true cones as ``large'' sets. Note that one must limit one's attention to \emph{true} cones for this to make sense. Indeed, $[\bot]_T$ is a cone but it is not large since it is a \emph{single} Lindenbaum degree.

We must introduce one more technical condition before stating the theorem.
\begin{definition}
A function $\mathfrak{g}$ is \emph{bounded} if there exists an $n\in\mathbb{N}$ such that for all $\varphi$, $\mathfrak{g}(\varphi)$ is $\Pi_n$.
\end{definition}

\begin{theorem}[\cite{walsh2020note}]\label{main}
Let $\mathfrak{g}$ be recursive, monotone, and bounded. Then one of the following holds:
\begin{itemize}
\item There is a true cone $\mathfrak{C}$ such that for all $\varphi \in \mathfrak{C}$, $T+\varphi\vdash\mathfrak{g}(\varphi).$
\item There is a true cone $\mathfrak{C}$ such that for all $\varphi \in \mathfrak{C}$, $T+\varphi+\mathfrak{g}(\varphi) \vdash \mathsf{Con}_T(\varphi).$
\end{itemize}
\end{theorem}
Roughly, this result states that any bounded recursive monotone operator must either be as weak as the identity operator in the limit or as strong as the consistency operator in the limit. By analogy with Sacks' question, this result confirms the suspicion that there are no natural theories of deductive strength strictly between $T$ and $T+\mathsf{Con}(T)$. If there were one, then relativizing it would produce a monotone operator of strictly intermediate strength, but there is no such operator.\footnote{In fact, this conclusion can be inferred from Theorem \ref{mw-first}, which was proved first.}

Logicians have studied many variants of the consistency operator that display interesting logical behavior. Sometimes these variants are mistaken for counter-examples to Theorem \ref{main}, so it is worth pausing to discuss why this is not the case.  

First, let's consider the case of slow reflection over $\mathsf{PA}$ \parencite{friedman2013slow}. The slow consistency statement for $\mathsf{PA}$ has the following form:
$$\mathsf{Con}^\diamondsuit( \mathsf{PA}) := \forall x \big( F_{\varepsilon_0}(x) \downarrow \to \mathsf{Con}({I\Sigma_x})\big).$$
Note that $F_{\varepsilon_0}$ is a recursive function whose totality is not provable in $\mathsf{PA}$. Since $F_{\varepsilon_0}$ is recursive, its graph is $\Sigma_1$-definable, whence $\mathsf{Con}^\diamondsuit(\mathsf{PA})$ is a $\Pi_1$ sentence. The slow consistency statement is not provable in $\mathsf{PA}$, but it is strictly weaker statement than the ordinary consistency statement $\mathsf{Con}(\mathsf{PA})$. That is:
\begin{theorem}[\cite{friedman2013slow}]
We have each of the following:
\begin{itemize}
        \item $\mathsf{PA}\vdash \mathsf{Con}(\mathsf{PA}) \to \mathsf{Con}^\diamondsuit(\mathsf{PA})$
        \item $\mathsf{PA}\nvdash \mathsf{Con}^\diamondsuit(\mathsf{PA}) \to \mathsf{Con}(\mathsf{PA})$
        \item $\mathsf{PA}\nvdash \mathsf{Con}^\diamondsuit (\mathsf{PA})$
    \end{itemize}
\end{theorem}
Relativizing the definition of slow consistency yields a slow consistency operator:
$$\mathsf{Con}^\diamondsuit_{\mathsf{PA}}(\ulcorner\varphi\urcorner) := \forall x \big( F_{\varepsilon_0}(x) \downarrow \to \mathsf{Con}_{I\Sigma_x}(\ulcorner\varphi\urcorner)\big).$$
So why doesn't the existence of the slow consistency operator refute Theorem \ref{main}? The key point is that the slow consistency operator ``catches up'' to the standard consistency operator for all inputs $\varphi$ such that $\mathsf{PA}+\varphi \vdash F_{\varepsilon_0}\downarrow$. That is, if $\mathsf{PA}+\varphi \vdash F_{\varepsilon_0}\downarrow$ then also $\mathsf{PA}\vdash \big(\varphi \wedge \mathsf{Con}_{\mathsf{PA}}(\ulcorner\varphi\urcorner)\big)\leftrightarrow \big(\varphi \wedge \mathsf{Con}^\diamondsuit_{\mathsf{PA}}(\ulcorner\varphi\urcorner)\big).$

Similar remarks explain why \emph{cut-free consistency} is not a counter-example over the base theory $\mathsf{EA}$. Cut-free consistency just means consistency with respect to proofs that do not use the cut rule. Of course, consistency and cut-free consistency coincide by the cut-elimination theorem. However, the cut-elimination algorithm is superexponential, and it cannot be proved to terminate in $\mathsf{EA}$. In fact, cut-free consistency exhibits the following curious behavior within $\mathsf{EA}$:
\begin{theorem}[\cite{visser1990interpretability}]
    For any sentence $\varphi$, 
    $$\mathsf{EA}\vdash \mathsf{Con}^{\textsc{cf}}_{\mathsf{EA}}\big( \mathsf{Con}^{\textsc{cf}}_{\mathsf{EA}}(\ulcorner\varphi\urcorner) \big)\leftrightarrow \mathsf{Con}_{\mathsf{EA}}(\ulcorner\varphi\urcorner).$$
\end{theorem}
Why isn't this a counter-example to Theorem \ref{main}? For any $\varphi$ which implies the cut-elimination theorem over $\mathsf{EA}$:
$$\mathsf{EA}\vdash  \big(\varphi \wedge \mathsf{Con}^{\textsc{cf}}_{\mathsf{EA}}(\ulcorner\varphi\urcorner)\big) \leftrightarrow \big(\varphi \wedge \mathsf{Con}_{\mathsf{EA}}(\ulcorner\varphi\urcorner)\big).$$ 

Slow consistency and cut-free consistency may seem to obviously be variants of the consistency operator. Indeed, both are defined by \emph{modifying} the definition of consistency (in one instance by ``slowing'' it down and in the other by removing a rule from the proof system). Theorem \ref{main} suggests one way in which this claim can be made precise. Slow consistency and cut-free consistency belong to the same equivalence class (modulo agreeing on a true cone) as the consistency operator.

Modulo an appropriate equivalence relation, Theorem \ref{main} illustrates that the consistency operator is canonical. It thereby indicates that there are no natural theories with deductive strength strictly between that of $T$ and $T+\mathsf{Con}_T$. Nevertheless, there are many ways one could imagine improving Theorem \ref{main}. For instance: Must $\mathfrak{g}$ be monotone, or would it be enough for $\mathfrak{g}$ to be extensional? Must $\mathfrak{g}$ be recursive, or would it be enough for $\mathfrak{g}$ to be limit-recursive? Finally, can the result be generalized to the iterates of the consistency operator? Let's look at each of these questions in turn.


\subsection{Extensional Operators}

One conspicuous feature of Theorem \ref{main} is that it concerns only monotone functions. Recall that Sacks' question and Martin's Conjecture concern degree-invariant functions, not order-preserving functions. Whereas monotone functions are analogues of order-preserving functions, extensional functions are analogues of degree-invariant functions. Indeed, extensional functions are merely functions that are invariant with respect to degrees in Lindenbaum algebras. So a strengthening of Theorem \ref{main} to extensional functions would yield a more faithful proof-theoretic analogue of Martin's Conjecture.

By a result of Shavrukov and Visser, Theorem \ref{main} cannot be strengthened to cover extensional functions in this way. S.\ Friedman asked them, in conversation, whether ``the consistency operator is, in some sense, the least way to strengthen theories.'' In response, Shavrukov and Visser exhibited an extensional density function on the Lindenbaum algebra of $\mathsf{PA}$. Their result actually extends to all reasonable theories by a result of Kripke and Pour-El \parencite{kripke1967deduction}.\footnote{Recall that a reasonable theory is a consistent, recursively axiomatized theory that interprets a modicum of arithmetic.}
\begin{theorem}[\cite{shavrukov2014uniform}]\label{shav-visser}
For any reasonable theory $T$, the Lindenbaum algebra of $T$ admits a recursive extensional density function.
\end{theorem}



Let me briefly explain the significance of this result. Rosser's theorem yields a recursive density function on the Lindenbaum algebra of $\mathsf{PA}$, yet this function is not extensional.\footnote{One can use Rosser-style self-reference to produce, for any sentences $\varphi$ and $\psi$ such that $\varphi$ strictly implies $\psi$, a sentence $\theta$ of strictly intermediate deductive strength. Since $\psi$ does not imply $\varphi$, this means that $\psi \wedge\neg \varphi$ is consistent. Hence, one can use Rosser-style self-reference to find a sentence $\theta$ independent from $\psi \wedge\neg \varphi$ and then take $\big(\varphi \vee (\psi \wedge \theta)\big)$ as one's intermediate sentence.} Here is how Shavrukov and Visser find an extensional density function: To start, they introduced the following formula (note that $\mathsf{Con}^2_T(\ulcorner\varphi\urcorner):=\mathsf{Con}\big(T+\varphi + \mathsf{Con}(T+\varphi)\big)$):
$$\mathsf{SV}(\varphi):= \varphi \wedge \forall x \big( \mathsf{Con}_{I\Sigma_x}(\ulcorner\varphi\urcorner) \to \mathsf{Con}^2_{I\Sigma_x}(\ulcorner\varphi\urcorner) \big)$$
The function $\varphi \mapsto \mathsf{SV}(\varphi)$ is extensional in the sense that whenever $\mathsf{PA}\vdash \varphi \leftrightarrow \psi$ then also $\mathsf{PA}\vdash \mathsf{SV}(\varphi) \leftrightarrow \mathsf{SV}(\psi)$. This formula yields a recursive extensional \emph{density} function on the Lindenbaum algebra of $\mathsf{PA}$. Indeed, we may define:
$$\mathsf{SV}^\star (\varphi,\psi) := \varphi \vee \big(\mathsf{SV}(\neg \varphi \wedge \psi) \wedge \psi \big).$$
This is an \emph{extensional} function in the sense that:
$$\mathsf{PA}\vdash (\varphi\leftrightarrow \psi) \wedge (\theta\leftrightarrow \zeta) \Rightarrow \mathsf{PA}\vdash \mathsf{SV}^\star(\varphi,\theta) \leftrightarrow \mathsf{SV}^\star(\psi,\zeta).$$
It is also a \emph{density} function in the sense that for any $\varphi$ and $\psi$ such that $\mathsf{PA}+\varphi\vdash \psi$ but $\mathsf{PA}+\psi\nvdash \varphi$:
\begin{enumerate}
    \item $\mathsf{PA}+ \varphi \vdash \mathsf{SV}^\star (\varphi,\psi)$
    \item $\mathsf{PA}+ \mathsf{SV}^\star (\varphi,\psi) \nvdash \varphi$
    \item $\mathsf{PA}+\mathsf{SV}^\star (\varphi,\psi) \vdash \psi$
    \item $\mathsf{PA}+\psi \nvdash \mathsf{SV}^\star (\varphi,\psi)$
\end{enumerate}

As a corollary we infer that the generalization of Theorem \ref{main} to extensional functions is false. Indeed, the function $\varphi \mapsto \mathsf{SV}^\star \big(\varphi \wedge\mathsf{Con}_\mathsf{PA}(\ulcorner\varphi\urcorner), \varphi\big)$ is everywhere \emph{stronger} than the identity but \emph{weaker} than the consistency operator. Thus, it is neither as weak as the identity operator nor as strong as the consistency operator on a true cone. In particular, Theorem \ref{shav-visser} shows that the function $\varphi \mapsto \mathsf{SV}^\star\big(\varphi \wedge\mathsf{Con}_\mathsf{PA}(\ulcorner\varphi\urcorner),\varphi\big)$ is an extensional function that is strictly stronger than the identity operator but strictly weaker than the consistency operator. Shavrukov and Visser wrote that this suggests a negative answer to S. Friedman's question whether the consistency operator is ``the least way to strengthen theories.''





\subsection{Recursiveness and Limit-recursiveness}

Must $\mathfrak{g}$ be recursive? Note that the restriction to recursive functions in the statement of Theorem \ref{main} is analogous to the restriction to $\mathsf{AD}$ in Martin's Conjecture. In both cases, these restrict our attention to ``simple'' functions. Note that if we dropped the restriction to recursive functions in the statement of Theorem \ref{main} we could easily find a counter-example using the Axiom of Choice. Yet there is a long distance in between \emph{recursive} functions and \emph{arbitrary} functions. Can Theorem \ref{main} be strengthened by assuming only that $\mathfrak{g}$ is arithmetic? Or by assuming only that $\mathfrak{g}$ is $0'$-recursive?

In \cite{walsh2020note} it is shown that the assumption that the function is recursive is necessary in the statement of Theorem \ref{main}. We do this by exhibiting a function $f$ that is limit-recursive but not recursive which meets the other hypotheses of Theorem \ref{main} but does not satisfy the conclusion. In particular, $f$ oscillates between behaving like the identity operator and the consistency operator, without converging on either. To make this precise, let's introduce some terminology.

\begin{definition}
Given a set $\mathfrak{A}$, we say that \emph{cofinally many true} sentences belong to $\mathfrak{A}$ if for every true $\varphi$ there is a true $\psi$ such that $T+\psi\vdash\varphi$ and $\psi \in \mathfrak{A}$.
\end{definition}
The motivation for the previous definition is this: If cofinally many true sentences belong to $\mathfrak{A}$, then the set of Lindenbaum degrees of members of $\mathfrak{A}$ is cofinal in the ultrafilter of Lindenbaum degrees of true sentences. 

\begin{theorem}[\cite{walsh2020note}]
There is a $0'$-recursive monotone $\mathfrak{g}$ such that, for all $\varphi$, $\mathfrak{g}(\varphi)$ is $\Pi_1$, and both of the following hold:
\begin{enumerate}
\item For cofinally many true $\varphi$: $[\varphi \wedge \mathfrak{g}(\varphi)]_T = [\varphi]_T.$
\item For cofinally many true $\varphi$: $[\varphi \wedge \mathfrak{g}(\varphi)]_T = [\varphi \wedge \mathsf{Con}_T(\ulcorner\varphi\urcorner)]_T.$
\end{enumerate}
\end{theorem}

By G\"odel's second incompleteness theorem, $[\varphi]_T\neq[\varphi \wedge \mathsf{Con}_T(\ulcorner\varphi\urcorner)]_T$ for $T$-consistent $\varphi$. So a function that oscillates between these two outputs cannot converge on either. That is, such a function can neither be as weak as the identity operator in the limit nor as strong as the consistency operator in the limit.



\subsection{Iterates of Consistency}

For the rest of this section we will consider the extent to which the features that make the consistency operator canonical generalize to iterates of the consistency operator. We define iterates of consistency so that they are always finite extensions of the base theory $T$ (that is, theories of the form $T+\varphi$ for some sentence $\varphi$). To define our iterates of consistency, we fix an effective presentation $\prec$ of a recursive ordinal. We informally define the iterates of $\mathsf{Con}_T$ so that for any sentence $\varphi$:
\begin{itemize}
    \item $\mathsf{Con}_T^0(\ulcorner\varphi\urcorner):=\top$
    \item $\mathsf{Con}_T^{\alpha+1}(\ulcorner\varphi\urcorner):=\mathsf{Con}_T(\ulcorner\varphi\wedge \mathsf{Con}_T^\alpha(\ulcorner\varphi\urcorner)\urcorner)$
    \item $\mathsf{Con}_T^\lambda(\ulcorner\varphi\urcorner):=\forall\alpha\prec\lambda \mathsf{Con}_T\big(\ulcorner\varphi \wedge \mathsf{Con}_T^\alpha(\ulcorner\varphi\urcorner)\urcorner\big)$
\end{itemize}
Strictly speaking, we define the iterates of $\mathsf{Con}_T$ via G\"{o}del's fixed point lemma:
$$T\vdash \mathsf{Con}_T^\alpha(\ulcorner\varphi\urcorner)\leftrightarrow\forall\beta\prec\alpha \; \mathsf{Con}_T\big(\ulcorner\varphi\wedge \mathsf{Con}_T^\beta(\ulcorner\varphi\urcorner)\urcorner\big)$$


\begin{remark}
Let us reiterate that all of the iterates of the consistency operator we have defined are sentences. This is appropriate for the current investigation since we are interested in functions that induce functions on the Lindenbaum algebra of $T$. For instance, given our definition, for a limit $\lambda$, $\mathsf{Con}^\lambda_T(\ulcorner\varphi\urcorner)$ is the sentence $\forall \alpha \prec \lambda \mathsf{Con}_T\big(\ulcorner\varphi \wedge \mathsf{Con}_T^\alpha(\ulcorner\varphi\urcorner)\urcorner\big)$. This clashes with the conventions adopted in some other papers, wherein $\mathsf{Con}^\lambda_T(\ulcorner\varphi\urcorner)$ is used as a name for the infinite set $\{ \mathsf{Con}^\alpha_T(\ulcorner\varphi\urcorner) \mid \alpha \prec \lambda \}$. For the same reason, it also comes apart from Turing's definition of theories axiomatized by iterated consistency statements, as introduced in \textsection \ref{turing-progressions}.
\end{remark}

Does Theorem \ref{main} generalize to iterates of the consistency operator along suitable ordinal notations?\footnote{We need only make some very meager assumptions on ordinal notations for the proofs of the results we cite to go through. Discussing these assumptions would be a distraction from the points we are trying to discuss here. For a definition of \emph{suitable} ordinal notations, see \cite[\textsection 2.2]{walsh2023evitable}.} An optimist might hope to prove the following:\\

\noindent \textbf{Na\"ive Conjecture:} Fix a suitable ordinal notation system $\prec$. Let $\alpha \succ 0$. Let $\mathfrak{g}$ be recursive and monotone such that, for all $\varphi$, $\mathfrak{g}(\varphi)$ is $\Pi_1$. Then one of the following holds:
\begin{enumerate}
\item There is a true cone $\mathfrak{C}$ such that for all $\varphi \in \mathfrak{C}$, $T+\varphi+\mathfrak{g}(\varphi) \vdash \mathsf{Con}^\alpha_T(\ulcorner\varphi\urcorner).$
\item For some $\beta\prec \alpha$, there is a true cone $\mathfrak{C}$ such that for all $\varphi \in \mathfrak{C}$,\\ $T+\varphi+\mathsf{Con}^\beta_T(\ulcorner\varphi\urcorner)\vdash\mathfrak{g}(\varphi).$
\end{enumerate}

As the conjecture's label suggests, the answer is negative. In fact, the proposed classification fails already at the very next step.

\begin{theorem}[\cite{walsh2023evitable}]\label{oscillating}
Fix a suitable ordinal notation system $\prec$. For every $\alpha\succ 0$, there is a recursive monotone $\mathfrak{g}$ such that, for all $\varphi$, $\mathfrak{g}(\varphi)$ is $\Pi_1$, and both of the following hold:
\begin{enumerate}
\item For cofinally many true $\varphi$: $[\varphi \wedge \mathfrak{g}(\varphi)]_T = [\varphi \wedge \mathsf{Con}^\alpha_T(\ulcorner\varphi\urcorner)]_T.$
\item For cofinally many true $\varphi$: $[\varphi \wedge \mathfrak{g}(\varphi)]_T=[\varphi \wedge\mathsf{Con}_T(\ulcorner\varphi\urcorner)]_T.$
\end{enumerate}
\end{theorem}

This is a dramatic failure of the na\"ive conjecture. Indeed, for $\alpha\succ 1$ and true $\varphi$, $[\varphi \wedge \mathsf{Con}^\alpha_T(\ulcorner\varphi\urcorner)]_T \neq [\varphi \wedge\mathsf{Con}_T(\ulcorner\varphi\urcorner)]_T$ by G\"odel's second incompleteness theorem. So a function that oscillates between mimicking consistency and mimicking $\alpha$-iterated consistency will neither be as weak as the consistency operator on a true cone nor as strong as any non-trivial iterate of the consistency operator on a true cone.


\subsection{More Positive Theorems}

Theorem \ref{oscillating} rules out a straightforward generalization of Theorem \ref{main} to the iterates of the consistency operator. Yet the counter-examples it produces still coincide with the iterates of the consistency operator on a non-trivial set of values. This is no coincidence. Indeed, every function that is, e.g., stronger than the consistency operator but as weak as the 2-iterated consistency operator must coincide with the latter on cofinally many true inputs. This latter claim generalizes upwards to all ordinal notations.

\begin{theorem}[\cite{montalban2019inevitability}]\label{mw-first}
Let $\mathfrak{g}$ be recursive and monotone. Let $\prec$ be a suitable ordinal notation system and let $\alpha$ be an ordinal notation from the field of $\prec$. Suppose that for all $\varphi$:
\begin{itemize}
\item $T+\varphi+ \mathsf{Con}^\alpha_T(\ulcorner\varphi\urcorner)\vdash \mathfrak{g}(\varphi)$. 
\item  For all $\beta\prec \alpha$:
\begin{itemize}
    \item $T+\varphi + \mathfrak{g}(\varphi)\vdash \mathsf{Con}^\beta(\ulcorner\varphi\urcorner)$.
    \item $T+\varphi + \mathsf{Con}^\beta(\ulcorner\varphi\urcorner) \nvdash \mathfrak{g}(\varphi) $.
\end{itemize}
\end{itemize}
Then for cofinally many true sentences $\varphi$:$$[\varphi \wedge \mathfrak{g}(\varphi)]_T = [\varphi\wedge \mathsf{Con}^\alpha_T(\ulcorner\varphi\urcorner)]_T.$$
\end{theorem}

Note that this rules out the possibility of recursive monotone functions with strength strictly between $\mathsf{Con}_T^\alpha$ and $\mathsf{Con}_T^{\alpha+1}$. Note, moreover, that Theorem \ref{mw-first} includes no restriction to bounded functions.

Theorem \ref{mw-first} says that if the range of a monotonic function $\mathfrak{g}$ is sufficiently constrained, then for
some $\varphi$ and some $\alpha$, 
$$[\mathfrak{g}(\varphi)]_T=[\varphi \wedge\mathsf{Con}^\alpha_T(\ulcorner\varphi\urcorner)]_T \neq [\bot]_T.$$
This property still holds even when these constraints on the range of $\mathfrak{g}$ are relaxed considerably. Indeed, the main theorem of \cite{montalban2019inevitability} states that if the strength of any sufficiently nice function $\mathfrak{g}$ is ``bounded'' by some iterate of the consistency operator, then $\mathfrak{g}$ must somewhere coincide with an iterate of the consistency operator.

\begin{theorem}[\cite{montalban2019inevitability}]\label{mw-second}
Let $\mathfrak{g}$ be recursive and monotone such that, for all $\varphi$, $\mathfrak{g}(\varphi)$ is $\Pi_1$. Let $\prec$ be a suitable ordinal notation system and let $\alpha$ be an ordinal notation from the field of $\prec$. 

Suppose that for all $\varphi$, $T+\varphi+\mathsf{Con}^\alpha_T(\ulcorner\varphi\urcorner)\vdash \mathfrak{g}(\varphi)$. Then, for some $\beta\preceq \alpha$ and some $\varphi$:
$$[\varphi \wedge \mathfrak{g}(\varphi)]_T=[\varphi \wedge \mathsf{Con}_T^\beta(\ulcorner\varphi\urcorner)]_T \neq[\bot]_T.$$
\end{theorem}

Theorem \ref{mw-second} was inspired by a recursion-theoretic result.

\begin{theorem}[\cite{slaman1988definable}]\label{slaman-steel}
Assume $\mathsf{AD}$. Suppose $f:2^\omega\to 2^\omega$ is order-preserving and increasing a.e.. Then either:
\begin{enumerate}
   \item[(i)] $x^\alpha <_T f(x)$ on a cone.
    \item[(ii)] For some $\beta \leq \alpha$, $f(x)\equiv_T x^\beta$ on a cone.
\end{enumerate}
\end{theorem} 

Theorem \ref{slaman-steel} is much stronger than Theorem \ref{mw-second} insofar as it achieves coincidence \emph{on a cone} rather than mere coincidence \emph{at some non-trivial point}. Theorem \ref{oscillating} rules out the possibility of strengthening Theorem \ref{mw-second} to secure coincidence on a cone. Yet other questions remain. 

\begin{question}
    Let $\mathfrak{g}$ be recursive and monotone such that, for all $\varphi$, $\mathfrak{g}(\varphi)$ is $\Pi_1$. 
Suppose that for all $\varphi$, $T+\varphi+\mathsf{Con}^\alpha_T(\ulcorner\varphi\urcorner)\vdash \mathfrak{g}(\varphi)$. 
Is there some \emph{true} $\varphi$ such that
$$[\varphi \wedge \mathfrak{g}(\varphi)]_T=[\varphi \wedge \mathsf{Con}_T^\beta(\ulcorner\varphi\urcorner)]_T?$$
\end{question}

The problem of formulating analogues of these results for reflection principles that are not $\Pi_1$-axiomatizable (such as 1-consistency) also remains open.


\section{Reflection Principles and Ordinal Analysis}\label{reflection-and-ordinals}

In this section we turn from consistency to other uniform reflection principles. For definitions of these principles see \textsection \ref{reflection-definitions}.

Recall that when Beklemishev (see \textsection \ref{Ordinal analysis via iterated reflection}) calculates the proof-theoretic ordinal of $T$, he first axiomatizes the $\Pi_{n+1}$ consequences of $T$ in terms of iterated $n$-consistency statements over a weak base theory. Not all theories can be studied using such methods, since not all theories are $\Pi_{n+1}$ equivalent to iterations of $n$-consistency. Nevertheless, ``natural'' theories seem to be $\Pi_{n+1}$ equivalent to iterations of $n$-consistency along ``natural'' well-orderings. The problem of characterizing the ``natural'' well-orderings is just as difficult as the problem of characterizing the ``natural'' theories. So the situation at present is this: There are many specific examples of connections between natural theories and iterations of reflection principles, but the latter depend on choices of ordinal notation system, so it is not clear whether these are part of any general connection. To prove such a general connection, we would need to study iterations of reflection principles in a way that is independent of ordinal notation systems. This is the theme that we pursue in this section.

\subsection{Descending Sequences in First-order Arithmetic}\label{ds-section}

To motivate our approach to this question, let's begin by discussing the apparent \emph{well-foundedness} of natural theories by proof-theoretic strength. There is a contrast to note between ascending sequences of theories and descending sequences of theories. Sequences of consistent theories that are ascending in consistency strength are quite easy to define: just iterate a strength-increasing principle along a well-ordering. Descending sequences---by contrast---seem to require more ingenuity to define. Does this reflect some inherent complexity of the descending sequences of theories in proof-theoretic hierarchies? In this section we turn to results, inspired by Steel's Theorem \ref{steel-descending}, to the effect that such descending sequence \emph{must} be complex. Indeed, these results show that descending sequences in proof-theoretic hierarchies are either complex or contain elements that are unsound.

Gaifman asked whether there is a recursive sequence $\langle T_0,T_1,\dots\rangle$ of consistent r.e.\ extensions of $\mathsf{PA}$ such that for all $n$, $T_n\vdash \mathsf{Con}(T_{n+1})$. The answer is positive; however, there are none that are $\mathsf{PA}$-\emph{provably} descending.

\begin{theorem}[H. Friedman, Smory\'nski, Solovay]\label{hss-theorem}
Both of the following hold:
\begin{enumerate}
    \item There is a recursive sequence $(T_n)_{n<\omega}$ of consistent r.e.\ extensions of $\mathsf{PA}$ such that for each $n$, $T_n\vdash \mathsf{Con}(T_{n+1}).$
    \item There is no r.e.\ sequence $(T_n)_{n<\omega}$ of consistent r.e.\ extensions of $\mathsf{PA}$ such that $\mathsf{PA}\vdash \forall x \mathsf{Pr}_{T_x}\big(\ulcorner\mathsf{Con}(T_{x+1})\urcorner\big)$.\footnote{This result is discussed in \cite{lindstrom2017aspects, smorynski2012self}. For a proof see \cite{pakhomov2021reflection}.}
\end{enumerate}
\end{theorem}

The second part of Theorem \ref{hss-theorem} takes a step towards the goal of this subsection. Indeed, it shows that any descending sequences in the consistency strength hierarchy over $\mathsf{PA}$ must be complicated in the following sense: Either it is not r.e.\ or it cannot be proved within $\mathsf{PA}$ to be a descending sequence.

The first part of Theorem \ref{hss-theorem} was improved by Visser.
\begin{theorem}[\cite{visser1988descending}]
There is a recursive sequence $(T_n)_{n<\omega}$ of consistent r.e.\ extensions of $\mathsf{PA}$ such for each $n$, $T_n$ proves each instance of the uniform reflection principle for $T_{n+1}$.
\end{theorem}
Recall that our goal in this subsection is to argue that descending sequences in reflection hierarchies are either complex or contain elements that are unsound. Visser's result shows that---without some condition of uniform provability---no such results can be proved with respect to complexity alone. However, if we keep \emph{both} complexity and unsoundness in mind, we can rule out descending sequences in reflection hierarchies, as the following result demonstrates:

\begin{theorem}[\cite{pakhomov2021reflection}]
\label{sigma2}
There is no recursively enumerable sequence $(T_n)_{n\in\mathbb{N}}$ of $\Sigma_2$-sound extensions of $B\Sigma_1$ such that, for each $n$, $T_n$ proves the $\Sigma_2$-soundness of $T_{n+1}$.\footnote{In \cite{pakhomov2021reflection} this is stated as a result about $\mathsf{EA}$ rather than $B\Sigma_1$. This is because we formalize our results in terms of ``smooth provability'' instead of using the ordinary provability predicate. Our proof works for smooth provability, but not ordinary provability, over $\mathsf{EA}$, whereas it works for both over $B\Sigma_1$. Rather than explicate the notion of smooth provability here we state our result about extensions of $B\Sigma_1$.}
\end{theorem}


Theorem \ref{sigma2} was inspired by an alternative proof, discovered by H. Friedman, of Steel's Theorem \ref{steel-descending} \parencite{steel1975descending, friedman1976uniformly}. For more recursion-theoretic applications of this technique see \cite{lutz2020incompleteness}. The adaptability of this proof arguably strengthens the analogy between the proof-theoretic and recursion-theoretic well-foundedness phenomena. 

There are a few ways that one can imagine refining Theorem \ref{sigma2}. In the first place, Theorem \ref{hss-theorem} concerns \emph{consistency} but Theorem \ref{sigma2} concerns $\Sigma_2$-soundness. These results skip right over $\Sigma_1$-soundness. So one might wonder whether the case of $\Sigma_1$-soundness is more like Theorem \ref{hss-theorem} or Theorem \ref{sigma2}. In fact, the situation with $\Sigma_1$-soundness is just like the situation with consistency.
\begin{theorem}[\cite{pakhomov2021reflection}]
\label{sigma1}
There exists a recursive sequence $(T_n)_{n\in\mathbb{N}}$ of $\Sigma_1$-sound extensions of elementary arithmetic such that, for each $n$, $T_n$ proves the $\Sigma_1$-soundness of $T_{n+1}$.
\end{theorem}
This demonstrates that Theorem \ref{sigma2} is the best possible theorem with respect to the syntactic complexity classes it concerns. 

In the second place, Theorem \ref{sigma2} is not clearly optimal with respect to its complexity assumptions. In particular, Theorem \ref{sigma2} shows that descending sequences in a certain proof-theoretic hierarchy cannot be \emph{simple} in the sense that they cannot be \emph{recursively enumerable}. However, the following question remains:
\begin{question}\label{ds-question}
Is there a limit-recursive sequence $(T_n)_{n\in\mathbb{N}}$ of $\Sigma_2$-sound extensions of $B\Sigma_1$ such that, for each $n$, $T_n$ proves the $\Sigma_2$-soundness of $T_{n+1}$?
\end{question}

\subsection{Descending Sequences in Second-order Arithmetic}

The method introduced to prove Theorem \ref{sigma2} easily adapts to rule out \emph{all} descending sequences of sufficiently sound theories in another hierarchy of proof-theoretic strength, namely, the one given by $\Pi^1_1$-reflection. Recall that a theory is $\Pi^1_1$-sound if all its $\Pi^1_1$ consequences are true.

\begin{theorem}[\cite{pakhomov2021reflection}]\label{pakhomov-walsh}
\label{well-foundedness_Pi11}
There is no sequence $(T_n)_{n\in\mathbb{N}}$ of $\Pi^1_1$-sound r.e.\ extensions of $\mathsf{ACA}_0$ such that, for each $n$, $T_n$ proves the $\Pi^1_1$-soundness of $T_{n+1}$.
\end{theorem}
 
Theorem \ref{well-foundedness_Pi11} rules out descending sequences of theories according to a certain metric of logical strength. It is worth pausing to note that when one proves that $U$ proves the consistency of $T$, one often also establishes the stronger fact that $U$ proves the $\Pi^1_1$-soundness of $T$. Consider, for instance, the construction of $\omega$-models of theories via forcing. Moreover, of course, many natural theories are $\Pi^1_1$-sound extensions of $\mathsf{ACA}_0$. Accordingly, Theorem \ref{well-foundedness_Pi11} delivers some insight into the phenomenon that the natural axiomatic theories are pre-well-ordered by consistency strength.

The following result provides an alternate perspective on the well-foundedness phenomenon.
\begin{theorem}[\cite{walsh2022incompleteness}]\label{walsh-thm}
There is no sequence $(T_n)_{n<\omega}$ of $\Pi^1_1$-sound and $\Sigma^1_1$-definable extensions of $\Sigma^1_1\text{-}\mathsf{AC}_0$ such that for each $n$, $T_n \vdash \mathsf{RFN}_{\Pi^1_1}(T_{n+1})$.
\end{theorem}
Note that whereas Theorem \ref{pakhomov-walsh} concerns r.e.\ theories, Theorem \ref{walsh-thm} concerns all $\Sigma^1_1$-definable theories. On the other hand, Theorem \ref{pakhomov-walsh} covers extensions of $\mathsf{ACA}_0$ whereas Theorem \ref{walsh-thm} applies only to extensions of the stronger system $\Sigma^1_1\text{-}\mathsf{AC}_0$. So Theorem \ref{walsh-thm} is neither stronger nor weaker than Theorem \ref{pakhomov-walsh} but rather incomparable with it.

It is notable that Theorem \ref{pakhomov-walsh} and Theorem \ref{walsh-thm} have totally different proofs. Theorem \ref{pakhomov-walsh} is proved by appealing to G\"odel's second incompleteness, but Theorem \ref{walsh-thm} is proved using concepts from ordinal analysis. One additional benefit of Theorem \ref{walsh-thm} is that it yields, as a corollary, a version of G\"{o}del's second incompleteness theorem.
\begin{theorem}[\cite{walsh2022incompleteness}]\label{walsh-inc-thm}
There is no $\Pi^1_1$-sound $\Sigma^1_1$-definable extension $T$ of $\Sigma^1_1\text{-}\mathsf{AC}_0$ such that $T \vdash \mathsf{RFN}_{\Pi^1_1}(T)$.
\end{theorem}
For comparison, G\"{o}del's original theorem states that no sufficiently strong $\Pi^0_1$-sound $\Sigma^0_1$-definable theory $T$ proves $\mathsf{RFN}_{\Pi^0_1}(T)$.\footnote{Usually, G\"{o}del's theorem is stated in terms of recursive axiomatizability and consistency. However, recursive axiomatizability is equivalent to $\Sigma^0_1$-definability by Craig's Trick. Moreover, consistency is provably equivalent (in elementary arithmetic) to $\Pi^0_1$-soundness.} 


\subsection{Ordinal Analysis}

Recall that in ordinal analysis, values called \emph{proof-theoretic ordinals} are systematically assigned to theories. According to conventional wisdom, proof-theoretic ordinals measure the proof-theoretic strength of theories. We will now turn to results that vindicate this conventional wisdom. Indeed, Theorem \ref{pakhomov-walsh} makes it possible to rank theories according to their strength in terms of provable $\Pi^1_1$-reflection. This leads to the following definition:
\begin{definition}
    The \emph{reflection rank} of $T$ is the rank of $T$ in the restriction of the relation $\{ (T,U) : U\vdash \mathsf{RFN}_{\Pi^1_1}(T) \}$ to $\Pi^1_1$-sound extensions of $\mathsf{ACA}_0$.\footnote{Theories must really be understood intensionally here for $\mathsf{RFN}_{\Pi^1_1}(T)$ to be well-defined. So think of theories here as primitive recursive presentations of sets of axioms.}
\end{definition}

For most theories $T$, this ``reflection rank'' of $T$ equals the proof-theoretic ordinal of $T$. We make the notion of ``most theories'' precise---just as in Martin's Conjecture---by restricting our attention to an appropriate ``cone'' of theories. In particular, we will focus on extensions of $\mathsf{ACA}_0^+$, which is itself an axiomatic theory extending $\mathsf{ACA}_0$ with the axiom ``every set is contained in an $\omega$-model of $\mathsf{ACA}_0$.''


\begin{theorem}[\cite{pakhomov2021reflection, pakhomov2021infinitary}]\label{ordinals_on_a_cone}
For any $\Pi^1_1$-sound extension $T$ of $\mathsf{ACA}_0^+$, the reflection rank of $T$ equals the $\Pi^1_1$ proof-theoretic ordinal of $T$.\footnote{Further connections between $\Pi^1_1$-reflection and ordinal analysis are forged in \cite[\textsection 6]{pakhomov2023reducing}.}
\end{theorem}

Theorem \ref{ordinals_on_a_cone} makes contact with the well-foundedness phenomenon but not the pre-linearity phenomenon. Now we will turn to results that make contact with both of these phenomena.

\subsubsection{Two Proof-theoretic Orderings}

In the introduction, we introduced the consistency strength ordering modulo a fixed theory $B$. In different contexts, logicians choose different theories for the base theory $B$. Common choices include $\mathsf{EA}$, $\mathsf{PRA}$, $\mathsf{PA}$, and $\mathsf{ZFC}$. In this section, we have to choose a specific base theory. Note that since $\mathsf{ACA}_0$ is conservative over $\mathsf{PA}$, consistency strength over $\mathsf{ACA}_0$ and consistency strength over $\mathsf{PA}$ are actually the same. Thus, we refine our definition of consistency strength as follows: 
\begin{definition}
$T\leq_{\mathsf{Con}} U \Leftrightarrow \mathsf{ACA}_0 \vdash \mathsf{Con}(U) \to \mathsf{Con}(T).$ 
\end{definition}

Another common way to compare the strength of theories is to compare their $\Pi^0_1$ consequences. 

\begin{definition}
$T\subseteq_{\Pi^0_1} U \Leftrightarrow$  For every $\varphi\in\Pi^0_1$, if $T\vdash \varphi$ then $U\vdash \varphi$.
\end{definition}

It is often claimed that, when we restrict our attention to ``natural'' theories, the $\subseteq_{\Pi^0_1}$ ordering coincides with relative consistency strength \parencite{steel2014godel}. However, these notions do not coincide in general. For a counter-example, consider any consistent $T$ and let $R_T$ be the Rosser sentence for $T$. Then $T+R_T\not\subseteq_{\Pi^0_1} T$ but $\mathsf{ACA}_0\vdash \mathsf{Con}(T) \to \mathsf{Con}(T+R_T)$, i.e., $T+R_T \leq_{\mathsf{Con}}T$.

The orderings $\leq_{\mathsf{Con}}$ and $\subseteq_{\Pi^0_1}$ are neither pre-linear nor pre-well-founded. That is, in both orderings, there are incomparable elements and infinite descending sequences. Remarkably, both of these features disappear when we restrict our attention to the natural theories. The restriction of $\leq_{\mathsf{Con}}$ to natural theories coincides with the restriction of $\subseteq_{\Pi^0_1}$ to natural theories, and these restrictions are pre-well-orderings.

Does this phenomenon we observe when we restrict to natural theories reflect some mathematical reality or is it an illusion? We will present results in the rest of this section that speak in favor of the former answer. Recall that in \textsection \ref{turing-progressions}--\ref{Ordinal analysis via iterated reflection} we discussed the possibility that the natural theories can be characterized in terms of the natural well-orderings.\footnote{The thought is that natural theories are $\Pi_1$ equivalent to iterations of consistency along natural well-orderings. Likewise, one might conjecture that the natural theories are $\Pi_{n+1}$ equivalent to iterations of $n$-consistency along natural well-orderings and $\Pi^1_1$-equivalent to iterations of $\Pi^1_1$-soundness along natural well-orderings.} Yet we do not have a convincing mathematical definition of the ``natural'' well-orderings. In the rest of this section we will introduce analogues of the aforementioned consistency strength and $\Pi^0_1$ inclusion orderings. These analogues move us into the a context that abstracts away from the problem of pathological ordinal notation systems. These analogues of the classical notions of proof-theoretic strength \emph{actually} coincide and \emph{actually} pre-well-order theories. This suggests that the problem of pathological ordinal notation systems truly is the impediment to pre-well-ordering in the original context.

\subsubsection{An Analogue of $\subseteq_{\Pi^0_1}$}
First, let's introduce our analogue of $\subseteq_{\Pi^0_1}$. The key idea is that we replace the notion of \emph{provability} with the notion of \emph{provability in the presence of an oracle for $\Gamma$ truth.} A theory $T$ will prove a sentence $\varphi$ in the presence of such an oracle if $T+\psi$ proves $\varphi$ for some true $\Gamma$ sentence $\psi$. We introduce the following notation to capture this idea:

\begin{definition}
We write $T\vdash^{\Gamma}\varphi$ if there is a true $\psi\in\Gamma$ such that $T+\psi\vdash\varphi$.
\end{definition}
Just as we focus on provability in the presence of an oracle, we focus not on \emph{inclusion of $\Gamma$ theorems} but on \emph{inclusion of $\Gamma$ theorems in the presence of an oracle}. We introduce the following notation for this notion:
\begin{definition}
We write $T\subseteq^{\Gamma}_{\Delta} U$ if for all $\varphi \in \Delta$, if $T\vdash ^{\Gamma} \varphi$ then  $U\vdash^{\Gamma}  \varphi$.
\end{definition}
The only theories we will consider are extensions of $\mathsf{ACA}_0$. Since $\mathsf{ACA}_0$ proves every true $\Sigma^0_1$ sentence,\footnote{Indeed, comparably weak subsystems of $\mathsf{ACA}_0$ are also $\Sigma^0_1$-complete.} we have the following equivalence:
$$T\vdash\varphi \Leftrightarrow T\vdash^{\Sigma^0_1}\varphi.$$ 
Thus, for all the theories we consider, we have the following equivalence:
$$T\subseteq_{\Pi^0_1}U \Leftrightarrow T\subseteq^{\Sigma^0_1}_{\Pi^0_1}U.$$
The analogue of $\subseteq^{\Sigma^0_1}_{\Pi^0_1}$ that we are interested in is $\subseteq^{\Sigma^1_1}_{\Pi^1_1}$, where:
$$T\subseteq^{\Sigma^1_1}_{\Pi^1_1} U \Leftrightarrow \text{For all $\varphi \in \Pi^1_1$, if $T\vdash ^{\Sigma^1_1} \varphi$ then  $U\vdash^{\Sigma^1_1}  \varphi$}$$
Note that the latter results from the former merely by replacing $\Sigma^0_1$  with $\Sigma^1_1$ and $\Pi^0_1$ with $\Pi^1_1$.


\subsubsection{An Analogue of $\leq_{\mathsf{Con}}$}

Before stating the analogue of $\leq_{\mathsf{Con}}$, let's note that (provably in $\mathsf{ACA}_0$) a theory is consistent just in case all of its $\Pi^0_1$ consequences are true.\footnote{Indeed, this is true of weak subsystems of $\mathsf{ACA}_0$ such as $\mathsf{EA}$.} Accordingly, we have the following equivalence:
$$T\leq_{\mathsf{Con}}U \Leftrightarrow \mathsf{ACA}_0 \vdash^{\Sigma^0_1} \mathsf{RFN}_{\Pi^0_1}(U) \to \mathsf{RFN}_{\Pi^0_1}(T),$$
where $\mathsf{RFN}_{\Pi^0_1}(T)$ is a formula expressing that all of $T$'s $\Pi^0_1$ consequences are true.

We thus generalize the notion of consistency strength as follows: 
\begin{definition}
Fixing a formula $\mathsf{RFN}_\Delta(T)$ expressing $\Delta$-soundness of $T$ in $\mathsf{ACA}_0$:
$$T\leq^{\Gamma}_{\mathsf{RFN_{\Delta}}} U \Leftrightarrow \mathsf{ACA}_0 \vdash^{\Gamma} \mathsf{RFN}_{\Delta}(U) \to \mathsf{RFN}_{\Delta}(T).$$
\end{definition}

We will be particularly interested in $\Pi^1_1$-soundness, where a theory is $\Pi^1_1$-sound just in case all its $\Pi^1_1$ consequences are true. We can formalize the $\Pi^1_1$-soundness of $T$ with a single sentence in $\mathsf{ACA}_0$:
$$\mathsf{RFN}_{\Pi^1_1}(T) : =\forall \varphi \in \Pi^1_1 \big( \mathsf{Pr}_T(\varphi) \to \mathsf{True}_{\Pi^1_1}(\varphi)  \big).$$
Note that $\mathsf{Pr}_T$ here picks out ordinary provability from $T$, not provability in the presence of an oracle. Hence, we have the following analogue of $\leq_{\mathsf{Con}}$:
$$T\leq^{\Sigma^1_1}_{\mathsf{RFN_{\Pi^1_1}}} U \Leftrightarrow \mathsf{ACA}_0 \vdash^{\Sigma^1_1} \mathsf{RFN}_{\Pi^1_1}(U) \to \mathsf{RFN}_{\Pi^1_1}(T).$$
Note once again that the definition of $\leq^{\Sigma^1_1}_{\mathsf{RFN_{\Pi^1_1}}}$ results from our equivalence characterizing $\leq_{\mathsf{Con}}$ merely by replacing $\Sigma^0_1$  with $\Sigma^1_1$ and $\Pi^0_1$ with $\Pi^1_1$.

\subsubsection{Connecting the Orderings via Ordinal Analysis}

Recall that, according to conventional wisdom, calculating the proof-theoretic ordinal of a theory is a means of measuring its logical strength. However, note that the ordering of theories induced by ordinal analysis:
$$T\leq_{\mathsf{WF} }U \Leftrightarrow |T|_{\mathsf{WF}}\leq |U|_{\mathsf{WF}}$$
is a pre-well-ordering since the ordinals are well-ordered. Hence, $\leq_{\mathsf{WF} }$ cannot strictly coincide with either $\leq_{\mathsf{Con}}$ or $\subseteq_{\Pi^0_1}$.

Nevertheless, in the presence of an oracle for $\Sigma^1_1$ truths, we can vindicate the common wisdom that ordinal analysis is a means of measuring the logical strength of theories:

\begin{theorem}[\cite{walsh2023characterizations}]\label{well-ordering}
For $\Pi^1_1$-sound arithmetically definable $T$ and $U$ extending $\mathsf{ACA}_0:$
$$T\subseteq_{\Pi^1_1}^{\Sigma^1_1}U \Longleftrightarrow |T|_{\mathsf{WF}}\leq|U|_{\mathsf{WF}}  \Longleftrightarrow T\leq^{\Sigma^1_1}_{\mathsf{RFN}_{\Pi^1_1}} U.$$
\end{theorem}
Theorem \ref{well-ordering} cannot be extended to $\Sigma^1_1$-definable theories. Since the ordinals are well-ordered, this immediately yields the following corollary:
\begin{corollary}
The relations $\subseteq_{\Pi^1_1}^{\Sigma^1_1}$ and $\leq^{\Sigma^1_1}_{\mathsf{RFN}_{\Pi^1_1}}$ pre-well-order the $\Pi^1_1$-sound arithmetically definable extensions of $\mathsf{ACA}_0$.
\end{corollary}
Note that in the statement of Theorem \ref{well-ordering} and its corollary, we have dropped the non-mathematical quantification over ``natural'' theories. \cite{jeon2023generalized} extends Theorem \ref{well-ordering} to set-theoretic reflection principles using a generalization of ordinal analysis due to \cite{pohlers1998subsystems}. 


\section{Conclusions}\label{conclusions}

Why are the natural axiomatic theories pre-well-ordered by proof-theoretic strength? We have described a few results that point towards an answer, interweaving the three themes discussed in \textsection \ref{interlude}. The emerging picture is that natural theories are proof-theoretically equivalent to iterated reflection principles. Moreover, ordinal analysis is a means of calculating the ranks of natural theories in these hierarchies. I write ``hierarchies'' since there will be different hierarchies attending different notions of proof-theoretic strength, depending, e.g., on the reflection principle they are based on. One may have to adapt the notion of ordinal analysis to investigate these different notions of proof-theoretic strength. The claim that ordinal analysis is a means of calculating the ranks of natural theories should be refined by further developing the systematic connections between iterated reflection and proof-theoretic ordinals.


\printbibliography

\end{document}